\documentclass{article}
\usepackage{amsmath,times}
\usepackage{amsthm}
\usepackage{verbatim}
\usepackage{amssymb}
\usepackage{latexsym}
\usepackage{graphicx}
\usepackage{amscd}
\usepackage[english]{babel}
\usepackage[latin1]{inputenc}
\newtheorem{theorem}{Theorem}[section]
\newtheorem{lemma}[theorem]{Lemma}
\newtheorem{proposition}[theorem]{Proposition}
\newtheorem{definition}[theorem]{Definition}
\newtheorem{corollary}[theorem]{Corollary}

\newtheorem{remark}[theorem]{Remark}
\begin{document} 

\title{Transforming metrics on a line bundle to the Okounkov body}

\date{\today{}}

\author{David Witt Nystr\"{o}m}

\maketitle

\begin{abstract}
Let $L$ be a big holomorphic line bundle on a complex projective manifold $X.$ We show how to associate a convex function on the Okounkov body of $L$ to any continuous metric $\psi$ on $L.$ We will call this the Chebyshev transform of $\psi,$ denoted by $c[\psi].$ Our main theorem states that the difference of metric volume of $L$ with respect to two metrics, a notion introduced by Berman-Boucksom, is equal to the integral over the Okounkov body of the difference of the Chebyshev transforms of the metrics. When the metrics have positive curvature the metric volume coincides with the Monge-Amp\`ere energy, which is a well-known functional in K\"ahler-Einstein geometry and Arakelov geometry. We show that this can be seen as a generalization of classical results on Chebyshev constants and the Legendre transform of invariant metrics on toric manifolds. As an application we prove the differentiability of the metric volume in the cone of big metrized $\mathbb{R}$-divisors. This generalizes the result of Boucksom-Favre-Jonsson on the differentiability of the ordinary volume of big $\mathbb{R}$-divisors and the result of Berman-Boucksom on the differentiability of the metric volume when the underlying line bundle is fixed.

\end{abstract}

$ $
$ $\\
Class. math.: 32W20, 32Q15, 32U20 

\tableofcontents

\section{Introduction}

In \cite{Khov2} (\cite{Khov3} is a published shortened version) and \cite{Lazarsfeld} Kaveh-Khovanskii and Lazarsfeld-Musta\c{t}\u{a} initiated a systematic study of Okounkov bodies of divisors and more generally of linear series. Our goal is to contribute with an analytic viewpoint.

It was Okounkov who in his papers \cite{Okounkov} and \cite{Okounkov2} introduced a way of associating a convex body in $\mathbb{R}^n$ to any ample divisor on a $n$-dimensional projective variety. This convex body, called the Okounkov body of the divisor and denoted by $\Delta(L)$, can then be studied using convex geometry. It was recognized in \cite{Lazarsfeld} that the construction works for arbitrary big divisors.

We will restrict ourselves to a complex projective manifold X, and instead of divisors we will for the most part use the language of holomorphic line bundles. Because of this, in the construction of the Okounkov body, we prefer choosing local holomorphic coordinates instead of the equivalent use of a flag of subvarieties (see \cite{Lazarsfeld}). We use additive notation for line bundles, i.e. we will write $kL$ instead of $L^{\otimes k}$ for the k:th tensor power of $L.$ We will also use the additive notation for metrics. If $h$ is a hermitian metric on a line bundle, we may write it as $h=e^{-\psi},$ and in this paper we will denote that metric by $\psi$. Thus if $\psi$ is a metric on $L,$ $k\psi$ is a metric on $kL.$ The pair $(L,\psi)$ of a line bundle $L$ with a continuous metric $\psi$ will be called a metrized line bundle.  

The main motivation for studying Okounkov bodies has been their connection to the volume function on divisors. Recall that the volume of a line bundle $L$ is defined as $$\textrm{vol}(L):=\limsup_{k \to \infty}\frac{n!}{k^n}\textrm{dim}(H^0(kL)).$$ A line bundle is said to be big if it has positive volume. From here on, all line bundles $L$ we consider will be assumed to be big. By Theorem A in \cite{Lazarsfeld}, for any big line bundle $L$ it holds that $$\textrm{vol}(L)=n!\textrm{vol}_{\mathbb{R}^n}(\Delta(L)).$$ 

We are interested in studying certain functionals on the space of metrics on $L$ that refine $\textrm{vol}(L)$.

The notion of a metric volume of a metrized line bundle $(L,\psi)$ was introduced by Berman-Boucksom in \cite{Berman}. Given a metric $\psi$ one has a natural norm on the the spaces of holomorphic sections $H^0(kL),$ namely the supremum norm $$||s||_{k\psi,\infty}:=\sup\{|s(x)|e^{-k\psi(x)/2}: x\in X\}.$$ Let $\mathcal{B}^{\infty}(k\psi)\subseteq H^0(kL)$ be the unit ball with respect to this norm.

$H^0(kL)$ is a vector space, thus given a basis we can calculate the volume of $\mathcal{B}^{\infty}(k\psi)$ with respect to the associated Lebesgue measure. This will depend on the choice of basis, but given a reference metric $\varphi$ one can compute the quotient  $$\frac{\textrm{vol}(\mathcal{B}^{\infty}(k\psi))}{\textrm{vol}(\mathcal{B}^{\infty}(k\varphi))}$$ and this quantity will be invariant under the change of basis. The $k$:th $\mathcal{L}$-bifunctional is defined as $$\mathcal{L}_k(\psi,\varphi):=\frac{n!}{2k^{n+1}}\log \left(\frac{\textrm{vol}(\mathcal{B}^{\infty}(k\psi))}{\textrm{vol}(\mathcal{B}^{\infty}(k\varphi))}\right).$$

The metric volume of a metrized line bundle $(L,\psi),$ denoted by $\textrm{vol}(L,\psi,\varphi),$ is defined as the limit 
\begin{equation} \label{introlimit}
\textrm{vol}(L,\psi,\varphi):=\lim_{k\to \infty}\mathcal{L}_k(\psi,\varphi).
\end{equation}

\begin{remark}
In \cite{Berman} this quantity is called the \emph{energy at equilibrium}, but we have in this paper chosen to call it the metric volume in order to accentuate the close relationship with the ordinary volume of line bundles.
\end{remark}

The metric volume obviously depends on the choice of $\varphi$ as a reference metric but it is easy to see that the difference of metric volumes $\textrm{vol}(L,\psi,\varphi)-\textrm{vol}(L,\psi',\varphi)$ is independent of the choice of reference.

The definition of the metric volume is clearly reminiscent of the definition of the volume of a line bundles. 
In fact, one easily checks that when adding $1$ to the reference metric $\varphi,$ we have that $$\textrm{vol}(L,\varphi+1,\varphi)=\textrm{vol}(L).$$ From this it follows readily that the metric volume is zero whenever the line bundle fails to be big.

In \cite{Berman} Berman-Boucksom prove that the limit (\ref{introlimit}) exists. They do this by proving that it actually converges to a certain integral over the space $X$ involving mixed Monge-Amp\`ere measures related to the metrics.

A metric $\psi$ is said to be \emph{psh} if the corresponding function expressed in a trivialization of the bundle is plurisubharmonic, so that $$dd^c\psi\geq 0$$ as a current. Given two locally bounded psh metrics $\psi$ and $\varphi$ one defines $\mathcal{E}(\psi,\varphi)$ as $$\frac{1}{n+1}\sum_{j=0}^n \int _X (\psi-\varphi)(dd^c\psi)^{j}\wedge (dd^c\varphi)^{n-j},$$ which we will refer to as the Monge-Amp\`ere energy of $\psi$ and $\varphi.$ This bifunctional first appeared in the works of Mabuchi and Aubin in K\"ahler-Einstein geometry (see \cite{Berman} and references therein). 

If $\psi$ and $\varphi$ are continuous but not necessarily psh, we may still define a Monge-Amp\`ere energy, by first projecting them down to the space of psh metrics, $$P(\psi):=\sup\{\psi' : \psi'\leq \psi, \psi' \textrm{ psh}\},$$ and then integrating over the Zariski-open subset $\Omega$ where the projected metrics are locally bounded. We are therefore led to consider the composite functional $\mathcal{E}\circ P:$
\begin{equation} \label{integralsimi}
\mathcal{E}\circ P(\psi,\varphi):=\frac{1}{n+1}\sum_{j=0}^n \int _{\Omega} (P(\psi)-P(\varphi))(dd^cP(\psi))^{j}\wedge (dd^cP(\varphi))^{n-j}.
\end{equation}
The Monge-Amp\`ere energy can also be seen as a generalization of the volume since if we let $\psi$ be equal to $\varphi+1,$ from e.g. \cite{Berman} we have that $$\mathcal{E}\circ P(\psi,\varphi)=\int_{\Omega}(dd^cP(\varphi))^n=\textrm{vol}(L).$$ This is not a coincidence. In fact Berman-Boucksom prove that for any pair of continuous metrics $\psi$ and $\varphi$ on a big line bundle $L$ we have that $$\mathcal{E}\circ P(\psi,\varphi)=\textrm{vol}(L,\psi,\varphi).$$

In \cite{Favre} Boucksom-Favre-Jonsson proved that the volume function on the N\'eron-Severi space is $\mathcal{C}^1$ in the big cone. This result was later reproved in \cite{Lazarsfeld} by Lazarsfeld-Musta\c{t}\u{a} using Okounkov bodies. Berman-Boucksom proved in \cite{Berman} the differentiability of the metric volume when the line bundle is fixed. A natural question is what one can say about the regularity of the metric volume when the line bundle is allowed to vary as well. In this paper we approach this question by combining the pluripotential methods of Berman-Boucksom with Okounkov body techniques inspired by the work of Lazarsfeld-Musta\c{t}\u{a}.
   
Given a continuous metric $\psi,$ we will show how to construct an associated convex function on the interior of the Okounkov body of $L$ which we will call the Chebyshev transform of $\psi,$ denoted by $c[\psi].$ The construction can be seen to generalize both the Chebyshev constants in classical potential theory and the Legendre transform of convex functions (see subsections 9.2 and 9.3 respectively).

First we describe how to construct $\Delta(L).$ Choose a point $p\in X$ and local holomorphic coordinates $z_1,...,z_n$ centered at $p.$ Choose also a trivialization of $L$ around $p.$ With respect to this trivialization any holomorphic section $s\in H^0(L)$ can be written as a convergent power series in the coordinates $z_i,$ $$s=\sum_{\alpha}a_{\alpha}z^{\alpha}.$$ Consider the lexicographic order on $\mathbb{N}^n,$ and let $v(s)$ denote the smallest index $\alpha$ (i.e. with respect to the lexicographic order) such that $$a_{\alpha} \neq 0.$$ We let $v(H^0(L))$ denote the set $\{v(s) : s\in H^0(L), s\neq 0\},$ and finally let the Okounkov body of $L$, denoted by $\Delta(L),$ be defined as closed convex hull in $\mathbb{R}^n$ of the union $$\bigcup_{k\geq 1}\frac{1}{k}v(H^0(kL)).$$ Observe that the construction depends on the choice of $p$ and the holomorphic coordinates. For other choices, the Okounkov bodies will in general differ.
  
Now let $\psi$ be a continuous metric on $L.$ There are associated supremum norms on the spaces of sections $H^0(kL),$ $$||s||_{k\psi}^2:=\sup_{x\in X} \{|s(x)|^2e^{-k\psi(x)}\}.$$  If $v(s)=k\alpha$ for some section $s\in H^0(kL),$ we let $A_{\alpha,k}$ denote the affine space of sections in $H^0(kL)$ of the form $$z^{k\alpha}+ \textrm{ higher order terms}.$$ We define the discrete Chebyshev transform $F[\psi]$ on $\bigcup_{k\geq 1} v(H^0(kL))\times \{k\}$ as $$F[\psi](k\alpha,k):=\inf\{\ln ||s||_{k\psi}^2 : s\in A_{\alpha,k}\}.$$ 

\begin{theorem} \label{theoreminintro}
For any point $p\in \Delta(L)^{\circ}$ and any sequence $\alpha(k)\in \frac{1}{k}v(H^0(kL))$ converging to $p,$ the limit $$\lim_{k\to \infty} \frac{1}{k}F[\psi](k\alpha(k),k)$$ exists and only depends on $p.$ We may therefore define the Chebyshev transform of $\psi$ by letting $$c[\psi](p):=\lim_{k\to \infty} \frac{1}{k}F[\psi](k\alpha(k),k),$$ for any sequence $\alpha(k)$ converging to $p.$ 
\end{theorem}

The main observation underlying the proof is the fact that the discrete Chebyshev transforms are subadditive. Our proof is thus very much inspired by the work of Zaharjuta, who in \cite{Zaharjuta} used subadditive functions on $\mathbb{N}^n$ when studying directional Chebyshev constants, and also by the article \cite{Bloom} where Bloom-Levenberg extend Zaharjutas results to a more general metrized setting, but still in $\mathbb{C}^n$ (we show in section 7 how to recover the formula of Bloom-Levenberg from Theorem \ref{theoreminintro}). Another inspiration comes from the work of Rumely-Lau-Varley in Arakelov geometry (see below).

We prove a general statement concerning subadditive functions on subsemigroups of $\mathbb{N}^{d}$ that extend Zaharjuta's results. 

\begin{theorem} \label{thethe27}
Let $\Gamma\subseteq \mathbb{N}^{d}$ be a semigroup which generates $\mathbb{Z}^{d}$ as a group, and let $F$ be a subadditive function on $\Gamma$ which is locally bounded from below by some linear function. Then for any sequence $\alpha(k) \in \Gamma$ such that $|\alpha(k)|\to \infty$ and $\frac{\alpha(k)}{|\alpha(k)|}\to p\in \Sigma(\Gamma)^{\circ}$ ($\Sigma(\Gamma)$ denotes the convex cone generated by $\Gamma$) for some point $p$ in the interior of $\Sigma(\Gamma),$ the limit $$\lim_{k\to \infty}\frac{F(\alpha(k))}{|\alpha(k)|}$$ exists and only depends on $F$ and $p.$ Furthermore the function $$c[F](p):=\lim_{k\to \infty}\frac{F(\alpha(k))}{|\alpha(k)|}$$ thus defined on $\Sigma(\Gamma)^{\circ}\cap \Sigma^{\circ}$ is convex.
\end{theorem}

Theorem \ref{theoreminintro} will follow from Theorem \ref{thethe27}. 

Our main result on the Chebyshev transform is the following.

\begin{theorem} \label{mainintro}
Let $\psi$ and $\varphi$ be two continuous metrics on $L.$ Then it holds that
\begin{equation} \label{introformula}
\textrm{vol}(L,\psi,\varphi)=n!\int_{\Delta(L)^{\circ}}(c[\varphi]-c[\psi])d\lambda,
\end{equation}
where $d\lambda$ denotes the Lebesgue measure on $\Delta(L).$
\end{theorem}  

The proof of Theorem \ref{mainintro} relies on the fact that one can use certain $L^2$-norms related to the metric, called Bernstein-Markov norms, to compute the Chebyshev transform. With the help of these one can interpret the right-hand side in equation (\ref{introformula}) as a limit of Donaldson bifunctionals  closely related and asymptotically equal to the ones used in the definition of the metric volume. This gives a new proof of the fact that the limit (\ref{introlimit}) exists.

In the setting of Arakelov geometry one studies adelic metrized line bundles $\overline{\mathcal{L}}$, and there is a corresponding notion of metric volume called sectional capacity. The relationship between these concepts is described in \cite{Berman}. The sectional capacity is defined as a limit of volumes of adelic unit balls in the space of adelic sections of powers of $\overline{\mathcal{L}}$. The existence of the limit was proved for ample adelic line bundles by Rumely-Lau-Varley in \cite{Rumely}. The method is similar to ours in that it defines a Chebyshev type transform following Zaharjuta's construction of directional Chebyshev constants. In order to define the directional Chebyshev constants Rumely-Lau-Varley constructs an ordered basis for the ring of sections with good multiplicative properties similiar to those of the monomial basis. In this paper we use the fact that the local holomorphic coordinates used when defining the Okounkov body also gives rise to a natural system of affine spaces with good multiplicative properties, and that this allows us to define our Chebyshev constants. For more on the use of Okounkov bodies in arithmetic geometry see \cite{Chen, Yuan2, Yuan}.  
  
Because of the homogeneity of the Okounkov body, i.e. $$\Delta(kL)=k\Delta(L),$$ one may define the Okounkov body of an arbitrary $\mathbb{Q}$-divisor $D$ by letting $$\Delta(D):=\frac{1}{p}\Delta(pD),$$ for any integer $p$ clearing all denominators in $D.$ Theorem B in \cite{Lazarsfeld} states that one may in fact associate an Okounkov body to an arbitrary big $\mathbb{R}$-divisor, such that the Okounkov bodies are fibers of a closed convex cone in $\mathbb{R}^n \times N^1(X)_{\mathbb{R}}$, where $N^1(X)_{\mathbb{R}}$ denotes the N\'eron-Severi space of $\mathbb{R}$-divisors. We show that this can be done also on the level of Chebyshev transforms, i.e. there is a continuous and indeed convex extension of the Chebyshev transforms to the space of continuous metrics on big $\mathbb{R}$-divisors. This shows that the metric volume also has a continuous extension to this space. 

As an application, using the differentiability result of Berman-Boucksom and some pluripotential theory and combining it with the new Okounkov body machinery we prove that the metric volume is differentiable. 

\begin{theorem} \label{theorem3inintro}
The metric volume function is $\mathcal{C}^1$ on the open cone of big $\mathbb{R}$-divisors equipped with two continuous metrics.
\end{theorem}

\subsection{Organization}

In section 2 we start by defining the Okounkov body of a semigroup, and we recall a result on semigroups by Khovanskii that will be of great use later on.

Section 3 deals with subadditive functions on subsemigroups of $\mathbb{N}^{n+1}$ and contains the proof of Theorem \ref{thethe27}.   

The definition of the Okounkov body of a line bundle follows in section 4. 

In section 5 we define the discrete Chebyshev transform of a metric, and prove that this function has the properties needed for Thereom  \ref{thethe27} to be applicable. We thus prove Theorem \ref{theoreminintro}. 

The metric volume is defined in section 6. Here we also state our main theorem, Theorem \ref{mainintro}.

In section 7 we show how one can use Bernstein-Markov norms instead of supremum norms in the construction of the Chebyshev transform.

The proof of Theorem \ref{mainintro} follows in section 8.

The Monge-Amp\`ere energy of metrics is introduced in section 9. We state the result of Berman-Boucksom which says that the metric volume is equal to a certain Monge-Amp\`ere energy. 

Section 10 discusses previuos results. 

In subsection 10.1 we observe that if we in (\ref{introformula}) let $\psi$ be equal to $\varphi+1,$ then we recover Theorem A in \cite{Lazarsfeld}, i.e. that $$\textrm{vol}(L)=n!\textrm{vol}_{\mathbb{R}^n}(\Delta(L)).$$ 

In subsection 10.2 we move on to clarify the connection to the classical Chebyshev constants. We see that if we embed $\mathbb{C}$ into $\mathbb{P}^1$ and choose our metrics wisely then formula (\ref{introformula}) gives us the classical result in potential theory that the Chebyshev constant and transfinite diamter of a regular compact set in $\mathbb{C}$ coincides. See subsection 9.2 for definitions.

Subsection 10.3 studies the case of a toric manifold, with a torus invariant line bundle and invariant metrics. We calculate the Chebyshev transforms, and observe that for invariant metrics, the Chebyshev transform equals the Legendre transform of the metric seen as a function on $\mathbb{R}^n$.  

We show in section 11 that if the first holomorphic coordinate $z_1$ defines a smooth submanifold $Y$ not contained in the augmented base locus of $L$ then the Chebyshev transform will be bounded near the interior of the zero-fiber of $\Delta(L),$ denoted by $\Delta(L)_0.$ It follows that the transform can be continuously extended to that part. 

We also note that one can define another  Cheyshev transform $c_{X|Y},$ defined on the interior of the zero fiber, by looking at a restricted subadditive function. When the line bundle is ample we prove, using the Ohsawa-Takegoshi extension theorem, that 
\begin{equation} \label{intrododo}
\mathcal{E}_Y(P(\varphi)_{|Y},P(\psi)_{|Y})=(n-1)!\int_{\Delta(L)_0} (c[\psi]-c[\varphi])(0, \alpha)d\alpha.
\end{equation}

In section 12 we show how to translate the results of Bloom-Levenberg to our language of Chebyshev transforms. We reprove Theorem 2.9 in \cite{Bloom} using our Theorem \ref{mainintro}, equation (\ref{intrododo}) and a recursion formula from \cite{Berman}.

We show in section 13 how to construct a convex and therefore continuous extension of the Chebyshev transform to arbitrary big $\mathbb{R}$-divisors.

In section 14 we move on to prove Theorem \ref{theorem3inintro} concerning the differentiability of the metric volume.

\subsection{Acknowledgement}

First of all I would like to thank Robert Berman for proposing the problem to me. I wish to thank Bo Berndtsson and S\'ebastien Boucksom for their numerous comments and suggestions concerning this article. Also I am very grateful for the insightful critique given to me by the anonymous referee.
 
 \section{The Okounkov body of a semigroup}

Let $\Gamma\subseteq \mathbb{N}^{n+1}$ be a subsemigroup of $\mathbb{N}^{n+1}$. We denote by $\Sigma(\Gamma)\subseteq \mathbb{R}^{n+1}$ the closed convex cone spanned by $\Gamma.$ By $\Delta_k(\Gamma)$ we will denote the set $$\Delta_k(\Gamma):=\{\alpha: (k\alpha,k)\in \Gamma\}\subseteq \mathbb{R}^{n}.$$

\begin{definition}
The Okounkov body $\Delta(\Gamma)$ of the semigroup $\Gamma$ is defined as $$\Delta(\Gamma):=\{\alpha: (\alpha,1)\in \Sigma(\Gamma)\}\subseteq  \mathbb{R}^{n}.$$
\end{definition}

It is clear that for all non-negative $k,$ $$\Delta_k(\Gamma)\subseteq \Delta(\Gamma).$$ The next theorem is a result of Khovanskii from \cite{Khovanskii}.

\begin{theorem} \label{khov}
Assume that $\Gamma \subseteq \mathbb{N}^{n+1}$ is a finitely generated semigroup which generates $\mathbb{Z}^{n+1}$ as a group. Then there exists an element $z\in \Sigma(\Gamma),$ such that $$(z+\Sigma(\Gamma)) \cap \mathbb{Z}^{n+1} \subseteq \Gamma.$$
\end{theorem}

When working with Okounkov bodies of semigroups it is sometimes useful to reformulate Theorem \ref{khov} into the following lemma. 

\begin{lemma} \label{goodygoody}
Suppose that $\Gamma$ is finitely generated, generates $\mathbb{Z}^{n+1}$ as a group, and also that $\Delta(\Gamma)$ is bounded. Then there exists a constant $C$ such that for all $k,$ if $$\alpha \in \Delta(\Gamma)\cap \left(\frac{1}{k}\mathbb{Z}\right)^{n}$$ and if the distance between $\alpha$ and the boundary of $\Delta(\Gamma)$ is greater than $C/k,$ then in fact we have that $$\alpha \in \Delta_k(\Gamma).$$
\end{lemma}

\begin{proof}
By definition we that $$\alpha \in \Delta(\Gamma)\cap \left(\frac{1}{k}\mathbb{Z}\right)^{n} \qquad{} \textrm{iff} \qquad{} (k\alpha,k)\in \Sigma(\Gamma)\cap \mathbb{Z}^{n+1}.$$ Also by definition $$\alpha \in \Delta_k(\Gamma) \qquad{} \textrm{iff} \qquad{} (k\alpha,k)\in \Gamma.$$ By Theorem \ref{khov}, there exists $z\in \Sigma(\Gamma)$ such that $$(k\alpha,k)\in \Gamma \qquad{} \textrm{if} \qquad{} (k\alpha,k)-z\in \Sigma(\Gamma),$$ and since $\Sigma(\Gamma)$ is a cone, $(k\alpha,k)-z\in \Sigma(\Gamma)$ iff $(\alpha,1)-z/k\in \Sigma(\Gamma).$ If $(\alpha,1)$ lies further than $|z|/k$ from the boundary of $\Sigma(\Gamma),$ then trivially $(\alpha,1)-z/k\in \Sigma(\Gamma).$ Since by assumtion the Okounkov body is bounded, the distance between $(\alpha,1)$ and the boundary of $\Sigma(\Gamma)$ is greater than some constant times the distance between $\alpha$ and the boundary of $\Delta(\Gamma).$ The lemma follows.
\end{proof}

\begin{corollary} \label{lattice}
Suppose that $\Gamma$ generates $\mathbb{Z}^{n+1}$ as a group, and also that $\Delta(\Gamma)$ is bounded. Then $\Delta(\Gamma)$ is equal to the closure of the union $\cup_{k\geq 0}\Delta_k(\Gamma).$
\end{corollary}

\begin{proof}
That $$\overline{\cup_{k\geq 0}\Delta_k(\Gamma)}\subseteq \Delta(\Gamma)$$ is clear. For the opposite direction, we exhaust $\Delta(\Gamma)$ by Okounkov bodies of finitely generated subsemigroups of $\Gamma$. Therefore, without loss of generality we may assume that $\Gamma$ is finitely generated. We apply Lemma \ref{goodygoody} which says that all the $(\frac{1}{k}\mathbb{Z})^{n}$ lattice points in $\Delta(\Gamma)$ whose distance to the boundary of $\Delta(\Gamma)$ is greater that some constant depending on the element $z$ in Theorem \ref{khov}, divided by $k$, actually lie in $\Delta_k(\Gamma)$. The corollary follows. 
\end{proof}

\section{Subadditive functions on semigroups}

Let $\Gamma$ be a semigroup. A real-valued function $F$ on $\Gamma$ is said to be \emph{subadditive} if for all $\alpha, \beta \in \Gamma$ it holds that 
\begin{equation} \label{subadditivity}
F(\alpha+\beta)\leq F(\alpha)+F(\beta).
\end{equation}

If $S$ is a subset of $\Gamma$ we say that a function $F$ is subadditive on $S$ if whenever $\alpha, \beta$ and $\alpha+\beta$ lie in $S$ the inequality (\ref{subadditivity}) holds.

If $\alpha \in \mathbb{R}^{n+1},$ we denote the sum of its coordinates $\sum \alpha_i$ by $|\alpha|.$ We also let $\Sigma^0\subseteq \mathbb{R}^{n+1}$ denote the set $$\Sigma^0:=\{(\alpha_1,...,\alpha_{n+1}) : |\alpha|=1, \alpha_i > 0\}.$$

In \cite{Bloom} Bloom-Levenberg observe that one can extract from \cite{Zaharjuta} the following theorem on subadditive functions on $\mathbb{N}^{n+1}.$

\begin{theorem}
Let $F$ be a subadditive function on $\mathbb{N}^{n+1}$ which is bounded from below by some linear function. Then for any sequence $\alpha(k)\in \mathbb{N}^{n+1}$ such that $|\alpha(k)| \to \infty$ when $k$ tends to infinity and such that $$\alpha(k)/|\alpha(k)| \to \theta \in \Sigma^0,$$ it holds that the limit $$c[F](\theta):=\lim_{k\to \infty}\frac{F(\alpha(k))}{|\alpha(k)|}$$ exists and does only depend on $\theta.$ Furthermore, the function $c[F]$ thus defined is convex on $\Sigma^0.$
\end{theorem}

We will give a proof of this theorem which also shows that it holds locally, i.e. that $F$ does not need to be subadditive on the whole of $\mathbb{N}^{n+1}$ but only on some open convex cone and only for large $|\alpha|.$ Then Zaharjuta's theorem still holds for the part of $\Sigma^0$ lying in the open cone. We will divide the proof into a couple of lemmas.

\begin{lemma} \label{golemma}
Let $O$ be an open convex cone in $\mathbb{R}_+^{n+1}$ and let $F$ be a subadditive function on $(O\setminus B(0,M))\cap \mathbb{N}^{n+1},$ where $B(0,M)$ denotes the ball of radius $M$ centered at the origin, and $M$ is any positive number. Then for any closed convex cone $K\subseteq O$ there exists a constant $C_K$ such that $$F(\alpha)\leq C_K|\alpha|$$ on $(K\setminus B(0,M))\cap \mathbb{N}^{n+1}.$  
\end{lemma}   

\begin{proof}
Pick finitely many points in $(O\setminus B(0,M))\cap \mathbb{N}^{n+1}$ such that if we denote by $\Gamma$ the semigroup generated by the points, the convex cone $\Sigma(\Gamma)$ should contain $(K\setminus B(0,M))$ and the distance between the boundaries should be positive. The points should also generate $\mathbb{Z}^{n+1}$ as a group. Then from Theorem \ref{khov} it follows that there exists an $M'$ such that 
\begin{equation} \label{rarara}
(K\setminus B(0,M'))\cap \mathbb{N}^{n+1}\subseteq \Gamma.
\end{equation}
Let $\alpha_i$ denote the generators of $\Gamma$ we picked. The inclusion (\ref{rarara}) means that for all $\alpha\in (K\setminus B(0,M'))\cap \mathbb{N}^{n+1}$ there exist non-negative integers $a_i$ such that $$\alpha=\sum a_i\alpha_i.$$ By the subadditivity we therefore get that $$F(\alpha)\leq \sum a_iF(\alpha_i)\leq C\sum a_i\leq C|\alpha|.$$ Since only finitely many points in $(K\setminus B(0,M))\cap \mathbb{N}^{n+1}$ do not lie in $(K\setminus B(0,M'))\cap \mathbb{N}^{n+1}$ the lemma follows. 
\end{proof}

\begin{lemma} \label{gogolemma}
Let $O, K$ and $F$ be as in the statement of Lemma \ref{golemma}. Let $\alpha$ be a point in $(K^{\circ}\setminus B(0,M))\cap \mathbb{N}^{n+1},$ and let $\gamma(k)$ be a sequence in $(K\setminus B(0,M))\cap \mathbb{N}^{n+1}$ such that $$|\gamma(k)|\to \infty$$ when $k$ tends to infinity and that $$\frac{\gamma(k)}{|\gamma(k)|} \rightarrow p\in K^{\circ}$$ for some point $p$ in the interior of $K.$ Let $l$ be the ray starting in $\alpha/|\alpha|,$ going through $p,$ and let $q$ denote the first intersection of $l$ with the boundary of $K.$ Denote by $t$ the number such that $$p=t\frac{\alpha}{|\alpha|}+(1-t)q.$$ Then there exists a constant $C_K$ depending only of $F$ and $K$ such that $$\limsup_{k\to \infty}\frac{F(\gamma(k))}{|\gamma(k)|}\leq t\frac{F(\alpha)}{|\alpha|}+(1-t)C_K.$$
\end{lemma}

\begin{proof}
We can pick points $\beta_i$ in $(K\setminus B(0,M))\cap \mathbb{N}^{n+1}$ with $\beta_i/|\beta_i|$ lying arbitrarily close to $q,$ such that if $\Gamma$ denotes the semigroup generated by the points $\beta_i$ and $\alpha,$ $\Gamma$ generates $\mathbb{Z}^{n+1}$ as a group and $$p\in \Sigma(\Gamma)^{\circ}.$$ Therefore from Theorem \ref{khov} it follows that for large $k$ $\gamma(k)$ can be written $$\gamma(k)=a\alpha+\sum a_i \beta_i$$ for non-negative integers $a_i$ and $a.$ The subadditivity of $F$ gives us that $$F(\gamma(k))\leq aF(\alpha)+\sum a_iF(\beta_i)\leq aF(\alpha)+C_K\sum a_i|\beta_i|,$$ where we in the last inequality used Lemma \ref{golemma}. Dividing by $|\gamma(k)|$ we get $$\frac{F(\gamma(k))}{|\gamma(k)|}\leq \frac{a|\alpha|}{|\gamma(k)|}\frac{F(\alpha)}{|\alpha|}+C_K \sum \frac{a_i|\beta_i|}{|\gamma(k)|}.$$ Our claim is that $\frac{a|\alpha|}{|\gamma(k)|}$ will tend to $t$ and that $\sum \frac{a_i|\beta_i|}{|\gamma(k)|}$ will tend to $(1-t).$ Consider the equations $$\frac{\gamma(k)}{|\gamma(k)|}=\frac{a|\alpha|}{|\gamma(k)|}\frac{\alpha}{|\alpha|}+\sum \frac{a_i|\beta_i|}{|\gamma(k)|} \frac{\beta_i}{|\beta_i|}$$ and $$p=t\frac{\alpha}{|\alpha|}+(1-t)q.$$ Observe that $$t=\frac{|p-\frac{\alpha}{|\alpha|}|}{|q-\alpha|}.$$ If $|\frac{\gamma(k)}{|\gamma(k)|}-p|<\delta$ and $|\frac{\beta_i}{|\beta_i|}-q|<\delta$ for all $i,$ then we see that $$\frac{a|\alpha|}{|\gamma(k)|}\leq \frac{|p-\frac{\alpha}{|\alpha|}|+\delta}{|q-\frac{\alpha}{|\alpha|}|-\delta}\leq t+\varepsilon(\delta),$$ where $\varepsilon(\delta)$ goes to zero as $\delta$ goes to zero. Similarly we have that 
\begin{equation} \label{rerera}
\frac{a|\alpha|}{|\gamma(k)|}\geq \frac{|p-\frac{\alpha}{|\alpha|}|-\delta}{|q-\frac{\alpha}{|\alpha|}|+\delta}\geq t- \varepsilon'(\delta),
\end{equation}
where $\varepsilon'(\delta)$ goes to zero as $\delta$ goes to zero. Since $$\frac{a|\alpha|}{|\gamma(k)|}+\sum \frac{a_i|\beta_i|}{|\gamma(k)|}=1,$$ inequality (\ref{rerera}) implies that $$\sum \frac{a_i|\beta_i|}{|\gamma(k)|}\leq 1-t+\varepsilon'(\delta).$$ The lemma follows.
\end{proof}

\begin{corollary} \label{gocoroll}
Let $O$ and $F$ be as in the statement of Lemma \ref{golemma}. Then for any sequence $\alpha(k)$ in $O\cap \mathbb{Z}^{n+1}$ such that $|\alpha(k)| \to \infty$ when $k$ tends to infinity and such that $\alpha(k)/|\alpha(k)|$ converges to some point $p$ in $O$ the limit $$\lim_{k \to \infty} \frac{F(\alpha)}{|\alpha(k)|}$$ exists and only depends on $F$ and $p.$
\end{corollary}

\begin{proof}
Let $\alpha(k)$ and $\beta(k)$ be two such sequences converging to $p.$ Let $K\subseteq O$ be some closed cone such that  $p\in K^{\circ}.$ Let us as in Lemma \ref{gogolemma} write $$p=t_k\frac{\beta(k)}{|\beta(k)|}+(1-t_k)q_k.$$ For any $\varepsilon > 0,$ $t_k$ is greater than $1-\varepsilon$ when $k$ is large enough. By Lemma \ref{gogolemma} we have that for such $k$ $$\limsup_{m\to \infty} \frac{F(\alpha(m))}{|\alpha(m)|}\leq (1-t_k)\frac{F(\beta(k))}{|\beta(k)|}+\varepsilon C_K\leq \frac{F(\beta(k))}{|\beta(k)|}+\varepsilon C_K +\varepsilon C,$$ where $C$ comes from the lower bound $$\frac{F(\beta)}{|\beta|}\geq C$$ which holds for all $\beta$ by assumption. Since $\varepsilon$ tends to zero when $k$ gets large we have that $$\limsup_{k\to \infty} \frac{F(\alpha(k))}{|\alpha(k)|}\leq \liminf_{k\to \infty} \frac{F(\beta(k))}{|\beta(k)|}.$$ By letting $\alpha(k)=\beta(k)$ we get existence of the limit, and by symmetry the limit is unique. 
\end{proof}

\begin{proposition} \label{propgo}
The function $c[F]$ on $O\cap \Sigma^{\circ}$ defined by $$c[F](p):=\lim_{k\to \infty}\frac{F(\alpha(k))}{|\alpha(k)|}$$ for any sequence $\alpha(k)$ such that $|\alpha(k)|\to \infty$ and $\frac{\alpha(k)}{|\alpha(k)|}\to p,$ which is well-defined according to Corollary \ref{gocoroll}, is convex, and therefore continuous. 
\end{proposition}

\begin{proof}
First we wish to show that $c[F]$ is lower semicontinuous. Let $p$ be a point in $O\cap \Sigma^{\circ}$ and $q_n$ a sequence converging to $p.$ From Lemma \ref{gogolemma} it follows that $$c[F](p) \leq \liminf_{q_n \rightarrow p}c[F](q_n),$$ which is equivalent to lower semicontinuity.

Using this the lemma will follow if we show that for any two points $p$ and $q$ in $O\cap \Sigma^{\circ}$ it holds that
\begin{equation} \label{Equation1000}
 2c[F](\frac{p+q}{2}) \leq c[F](p)+c[F](q).
\end{equation}

Choose sequences $\alpha(k), \beta(k) \in O\cap \mathbb{N}^{n+1}$ such that $$\frac{\alpha(k)}{|\alpha(k)|} \rightarrow p, \qquad \frac{\beta(k)}{|\beta(k)|} \rightarrow q,$$ and for simplicity assume that $|\alpha(k)|=|\beta(k)|.$ Then $$\frac{\alpha(k)+\beta(k)}{|\alpha(k)+\beta(k)|} \rightarrow \frac{p+q}{2}.$$ Hence 
\begin{eqnarray*}
2c[F](\frac{p+q}{2})=\lim_{k \rightarrow \infty}\frac{F(\alpha(k)+\beta(k))}{|\alpha(k)|}\leq \lim_{k \rightarrow \infty}\frac{F(\alpha(k))}{|\alpha(k)|}+\lim_{k \rightarrow \infty}\frac{F(\beta(k))}{|\beta(k)|}=\\=c[F](p)+c[F](q).
\end{eqnarray*}
\end{proof}

Together with Theorem \ref{khov} these lemmas yield a general result for subadditive functions on subsemigroups of $\mathbb{N}^{n+1}.$

A function $F$ defined on a cone $O$ is said to be \emph{locally linearly bounded from below} if for each point $p\in O$ there exists an open subcone $O'\subseteq O$ containing $p$ and a linear function $\lambda$ on $O'$ such that $F\geq \lambda$ on $O'.$  

\begin{theorem} \label{thethe}
Let $\Gamma\subseteq \mathbb{N}^{n+1}$ be a semigroup which generates $\mathbb{Z}^{n+1}$ as a group, and let $F$ be a subadditive function on $\Gamma$ which is locally linearly bounded from below. Then for any sequence $\alpha(k) \in \Gamma$ such that $|\alpha(k)|\to \infty$ and $\frac{\alpha(k)}{|\alpha(k)|}\to p\in \Sigma(\Gamma)^{\circ}$ for some point $p$ in the interior of $\Sigma(\Gamma),$ the limit $$\lim_{k\to \infty}\frac{F(\alpha(k))}{|\alpha(k)|}$$ exists and only depends on $F$ and $p.$ Furthermore the function $$c[F](p):=\lim_{k\to \infty}\frac{F(\alpha(k))}{|\alpha(k)|}$$ thus defined on $\Sigma(\Gamma)^{\circ}\cap \Sigma^{\circ}$ is convex.
\end{theorem}

\begin{proof}
By Theorem \ref{khov} it follows that for any point $p\in \Sigma(\Gamma)^{\circ}$ there exists an open convex cone $O$ and a number $M$ such that $$(O \setminus B(0,M))\cap \mathbb{N}^{n+1}\subseteq \Gamma.$$ We can also choose $O$ such that $F$ is bounded from below by a linear function on $O.$ Therefore the theorem follows immediately from Corollary \ref{gocoroll} and Proposition \ref{propgo}.
\end{proof}

We will show how this theorem can be seen as the counterpart to Theorem \ref{khov} for subadditive functions.

\begin{definition}
Let $\Gamma$ be a subsemigroup of $\mathbb{N}^{n+1}$ and let $F$ be a subadditive function of $\Gamma$ which is locally linearly bounded from below. One defines the convex envelope of $F,$ denoted by $P(F),$ as the supremum of all linear functions on $\Sigma(\Gamma)^{\circ}$ dominated by $F,$ or which ammounts to the same thing, the supremum of all convex one-homogeneous functions on $\Sigma(\Gamma)^{\circ}$ dominated by $F.$
\end{definition}

\begin{theorem} \label{jajatheo}
If $\Gamma$ generates $\mathbb{Z}^{n+1}$ as a group, then for any subadditive function $F$ on $\Gamma$ which is locally linearly bounded from below it holds that $$F(\alpha)= P(F)(\alpha)+o(|\alpha|)$$ for $\alpha\in \Gamma\cap \Sigma(\Gamma)^{\circ}.$
\end{theorem}

\begin{proof}
That $$F(\alpha)\geq P(F)(\alpha)$$ follows from the definition. If we let $c[F]$ be defined on the whole of $\Sigma(\Gamma)^{\circ}$ by letting $$c[F](\alpha):=|\alpha|c[F](\frac{\alpha}{|\alpha|}),$$ it follows from Theorem \ref{thethe} that $c[F]$ will be convex and one-homogeneous. It will also be dominated by $F$ since by the subadditivity $$\frac{F(\alpha)}{|\alpha|}\geq \frac{F(k\alpha)}{|k\alpha|}$$ for all positive integers and therefore $$\frac{F(\alpha)}{|\alpha|}\geq \lim_{k\to \infty}\frac{F(k\alpha)}{|k\alpha|}=c[F](\frac{\alpha}{|\alpha|}).$$ It follows that $$P(F)\geq c[F].$$ For $\alpha\in \Gamma$ by definition we have that $$P(F)(\alpha)\leq \frac{F(k\alpha)}{k}$$ for all positive integers $k.$ At the same time $$c[F](\alpha)=\lim_{k\to \infty}\frac{F(k\alpha)}{k},$$ hence we get that $$P(F)(\alpha)\leq c[F](\alpha)$$ for $\alpha\in \Gamma$ Since both $P(F)$ and $c[F]$ are convex they are continuous, so by the homogeneity we get that $$P(F)\leq c[F]$$ on $\Sigma(\Gamma)^{\circ}$, and therefore $P(F)=c[F].$ The theorem now follows from Theorem \ref{thethe}.
\end{proof}

\section{The Okounkov body of a line bundle}

In this section, following Okounkov, we will show how to associate a semigroup to a line bundle.

\begin{definition}
An order $<$ on $\mathbb{N}^n$ is additive if $\alpha <\beta$ and $\alpha' <\beta'$ implies that $$\alpha+\alpha'<\beta+\beta'.$$
\end{definition}
One example of an additive order is the lexicographic order where  $$(\alpha_1,...,\alpha_n)<_{\textrm{lex}}(\beta_1,...,\beta_n)$$ iff there exists an index $j$ such that $\alpha_j<\beta_j$ and $\alpha_i=\beta_i$ for $i<j.$

Let $X$ be a compact projective complex manifold of dimension $n$, and $L$ a holomorphic line bundle, which we will assume to be big. Suppose we have chosen a point $p$ in $X,$ and local holomorphic coordinates $z_1,...,z_n$ around that point, and let $e_p\in H^0(U,L)$ be a local trivialization  of $L$ around $p.$ Any holomorphic section $s\in H^0(X,kL)$ has an unique represention as a convergent power series in the variables $z_i,$ $$\frac{s}{e_p^k}=\sum a_{\alpha}z^{\alpha},$$ which for convenience we will simply write as $$s=\sum a_{\alpha}z^{\alpha}.$$ We consider the lexicographic order on the multiindices $\alpha$, and let $v(s)$ denote the smallest index $\alpha$ such that $a_{\alpha} \neq 0.$

\begin{definition}
Let $\Gamma(L)$ denote the set $$\bigcup_{k\geq 0}\left(v(H^0(kL))\times\{k\}\right) \subseteq \mathbb{N}^{n+1}.$$ It is a semigroup, since for $s\in H^0(kL)$ and $t\in H^0(mL)$ 
\begin{equation} \label{additive}
v(st)=v(s)+v(t).
\end{equation}
The Okounkov body of $L$, denoted by $\Delta(L)$, is defined as the Okounkov body of the associated semigroup $\Gamma(L).$
\end{definition}

We write $\Delta_k(\Gamma(L))$ simply as $\Delta_k(L).$

Let us recall some basic facts on Okounkov bodies.

\begin{lemma} \label{points}
The number of points in $\Delta_k(L)$ is equal to the dimension of the vector space $H^0(kL).$
\end{lemma}

This is part of Lemma 1.4 in \cite{Lazarsfeld}.

\begin{lemma} \label{okounkbound}
The Okounkov body of a big line bundle is bounded, hence compact.
\end{lemma}

This is Lemma 1.11 in \cite{Lazarsfeld}.

\begin{lemma} \label{simplex}
If $L$ is a big line bundle, $\Gamma(L)$ generates $\mathbb{Z}^{n+1}$ as a group. In fact $\Gamma(L)$ contains a translated unit simplex.
\end{lemma}

It is proved as part of Lemma 2.2 in \cite{Lazarsfeld}.

\begin{remark} \label{order}
Note that the additivity of $v$ as seen in equation (\ref{additive}) only depends on the fact that the lexicographic order is additive. Therefore we could have used any total additive order on $\mathbb{N}^n$ to define a semigroup $\tilde{\Gamma}(L),$ and the associated Okounkov body $\tilde{\Delta}(L).$ We will only consider the case where the Okounkov body $\tilde{\Delta}(L)$ is bounded, and the semigroup $\tilde{\Gamma}(L)$ generates $\mathbb{N}^n$ as a group.    
\end{remark}

\begin{lemma} \label{hull}
For any closed set $K$ contained in the convex hull of $\Delta_M(L)$ for some $M,$ there exists a constant $C_K$ such that if $$\alpha \in K\cap (\frac{1}{k}\mathbb{Z})^{n}$$ and the distance between $\alpha$ and the boundary of $K$ is greater than $\frac{C_K}{k},$ then $\alpha\in \Delta_k(L).$
\end{lemma}

\begin{proof} \label{lat}
Let $\Gamma$ be the semigroup generated by the elements $(M\beta,M)$ where $\beta \in \Delta_M(L),$ and some unit simplex in $\Gamma(L).$ Applying Lemma \ref{goodygoody} gives the lemma.
\end{proof}

\begin{lemma} \label{lattice2}
If $K$ is relatively compact in the interior of $\Delta(L),$ there exists a number $M$ such that for $k>M$, $$\alpha \in K\cap (\frac{1}{k}\mathbb{Z})^n$$ implies that $\alpha \in \Delta_k(L).$
\end{lemma}

\begin{proof}
This is a consequence of Lemma \ref{hull} by choosing $M$ such that the distance between $K$ and the convex hull of $\Delta_M(L)$ is strictly positive, therefore greater than $\frac{C_K}{k}$ for large $k$.
\end{proof}
  
\section{The Chebyshev transform}

\begin{definition}
A continuous hermitian metric $h=e^{-\psi}$ on a line bundle $L$ is a continuous choice of scalar product on the complex line $L_p$ at each point $p$ on the manifold. If $f$ is a local frame for $L$ on $U_f$, then one writes $$|f|^2=h_f=e^{-\psi_f},$$ where $\psi_f$ is a continuous  function on $U_f$. In this paper we let $\psi$ denote the metric $h=e^{-\psi}.$ 
\end{definition} 

The pair $(L,\psi)$ of a line bundle $L$ together with a continuous metric $\psi$ will be called a metrized line bundle.

We will show how to a given metrized line bundle $(L,\psi)$ one associates a subadditive function on the semigroup $\Gamma(L).$

For all $(k\alpha,k) \in \Gamma(L),$ let us denote by $A_{\alpha,k}$ the affine space of sections in $H^0(kL)$ of the form $$z^{k\alpha}+ \textrm{higher order terms}.$$ Consider the supremum norm $||.||_{k\psi}$ on $H^0(kL)$ given by $$||s||^2_{k\psi}:=\sup_{x\in X}\{|s(x)|^2e^{-k\psi(x)}\}.$$ 

\begin{definition}
We define the discrete Chebyshev transform $F[\psi]$ on $\Gamma(L)$ by  $$F[\psi](k\alpha,k):=\inf\{\ln ||s||_{k\psi}^2: s\in A_{\alpha,k}\}.$$
\end{definition} 

A section $s$ in $A_{\alpha,k}$ which minimizes the supremum norm is called a Chebyshev section.

\begin{lemma} \label{subsub}
The function $F[\psi]$ is subadditive.
\end{lemma}

\begin{proof}
Let $(k\alpha,k)$ and $(l\beta,l)$ be two points in $\Gamma(L),$ and denote by $\gamma$ $$\gamma:=\frac{k\alpha+l\beta}{k+l}.$$ Thus we have that $$(k\alpha,k)+(l\beta,l)=((k+l)\gamma,k+l).$$ Let $s$ be some section in $A_{\alpha,k}$ and $s'$ some section in $A_{\beta,l}.$ Since 
\begin{eqnarray*}
ss'=(z^{k\alpha}+ \textrm{higher order terms})(z^{l\beta}+ \textrm{higher order terms})=\\=z^{(k+l)\gamma}+\textrm{higher order terms},
\end{eqnarray*}
we see that $ss'\in A_{\gamma,k+l}.$ We also note that the supremum of the product of two functions is less or equal to the product of the supremums, i.e. $$||ss'||^2_{(k+l)\psi}\leq ||s||^2_{k\psi}||s'||^2_{l\psi}.$$  It follows that $$\inf\{ ||s||^2_{k\psi}: s\in A_{\alpha,k}\}\inf\{ ||s'||^2_{l\psi}: s'\in A_{\beta,l}\}\leq \inf\{||t||^2_{\gamma,k+l}: t\in A_{\gamma,k+l}\},$$ which gives the lemma by taking the logarithm.
\end{proof}

\begin{lemma} \label{lowerbound}
There exists a constant $C$ such that for all $(k\alpha,k)\in \Gamma(L),$ $$F[\psi](k\alpha,k)\geq C|(k\alpha,k)|$$.
\end{lemma}

\begin{proof}
Let $r>0$ be such that the polydisc $D$ of radius $r$ centered at $p$ is fully contained in the coordinate chart of $z_1,...,z_n.$ We can also assume that our trivialization $e_p\in H^0(U,L)$ of $L$ is defined on $D,$ i.e. $D\subseteq U.$ Let $s$ be a section in $A_{\alpha,k}$ and let $$\tilde{s}:=\frac{s}{e_p^k}.$$ Denote by $\psi_p$ the trivialization of $\psi.$ Hence $$|s|^2e^{-k\psi}=|\tilde{s}|^2e^{-k\psi_p}.$$ Since $\psi_p$ is continuous, $$e^{-\psi_p} > A$$ on $D$ for some constant $A.$ This yields that
\begin{eqnarray*}
||s||^2 \geq \sup_{x\in D} \{|\tilde{s}(x)|^2e^{-k\psi_p(x)}\}\geq A^k \sup_{x\in D} \{|\tilde{s}(x)|^2\}.
\end{eqnarray*}
We claim that $$\sup_{x\in D} \{|\tilde{s}(x)|^2\}\geq r^{k|\alpha|}.$$ Observe that $$\sup_{z\in D} \{|z^{k\alpha}|^2\}= r^{k|\alpha|}.$$ One now shows that $$\sup_{z\in D} \{|z^{k\alpha}|^2\}\leq \sup_{z\in D} \{|z^{k\alpha}+ \textrm{higher order terms}|^2\}$$ by simply reducing it to the case of one variable where it is immediate. We get that $$||s||^2\geq A^kr^{k|\alpha|}$$ and hence $$F[\psi](k\alpha,k)\geq k\ln A+k|\alpha|\ln r\geq C(k+k|\alpha|),$$ if we choose $C$ to be less than both $\ln A$ and $\ln r.$
\end{proof}

\begin{definition}
We define the Chebyshev transform of $\psi,$ denoted by $c[\psi]$ as the convex envelope of $F[\psi]$ on $\Sigma(\Gamma)^{\circ}.$ It is convex and one-homogeneous. We will also identify it with its restriction to $\Delta(L)^{\circ},$ the interior of the Okounkov body of L. Recall that by definition $$\Delta(L):=\Sigma(L)\cap (\mathbb{R}^n\times \{1\}).$$ \end{definition}

\begin{proposition} \label{propidopp}
For any sequence $(k\alpha(k),k)$ in $\Gamma(L),$ $k\to \infty,$ such that $$\lim_{k\to \infty}\alpha(k)=p\in \Delta(L)^{\circ},$$ it holds that $$c[\psi](p)=\lim_{k\to \infty}\frac{1}{k}\ln ||t_{\alpha(k),k}||^2,$$ where $t_{\alpha(k),k}$ is a Chebyshev section in $A_{\alpha((k),k}.$
\end{proposition}

\begin{proof}
By Lemma \ref{subsub} and Lemma \ref{lowerbound} we can apply Theorem \ref{jajatheo} to the function $F[\psi]$ and get that 
\begin{eqnarray*}
c[\psi](p)=|(p,1)|c[\psi](\frac{(p,1)}{|(p,1)|})=|(p,1)|\lim_{k\to \infty}\frac{F[\psi](k\alpha,k)}{k|(\alpha(k),1)|}=\\=\lim_{k\to \infty}\frac{F[\psi](k\alpha,k)}{k}=\lim_{k\to \infty}\frac{1}{k}\ln ||t_{\alpha(k),k}||^2.
\end{eqnarray*}
\end{proof}

\begin{lemma} \label{constantgo}
Let $\psi$ be a continuous metric on $L$ and consider the continuous metric on $L$ given by $\psi+C$ for some constant $C.$ Then it holds that 
\begin{equation} \label{eqeqeq}
F[\psi+C](k\alpha,k)=F[\psi](k\alpha,k)-kC,
\end{equation}
and that $$c[\psi+C]=c[\psi]-C$$ on $\Delta(L)^{\circ}.$  
\end{lemma}

\begin{proof}
For any section $s\in H^0(kL)$ we have that $$||s||^2_{k(\psi+C)}=e^{-kC}||s||^2_{k\psi},$$ therefore $$\ln ||s||^2_{k(\psi+C)}=\ln ||s||^2_{k\psi}-kC.$$ The lemma thus follows from the definitions.
\end{proof}

\begin{lemma} \label{monotone2}
If $\psi$ and $\varphi$ are two continuous metrics such that  $$\psi \leq \varphi,$$ then $$F[\psi]\geq F[\varphi],$$ and also $$c[\psi]\geq c[\varphi].$$
\end{lemma}

\begin{proof}
Follows immediately from the definitions.
\end{proof}

\begin{proposition}
For any two continuous metrics on $L,$ $\psi$ and $\varphi,$ the difference of the Chebyshev transforms, $c[\psi]-c[\varphi],$ is continuous and bounded on $\Delta(L)^{\circ}.$ 
\end{proposition}

\begin{proof}
It is the difference of two convex hence continuous functions, and is therefore continuous. Since $\psi-\varphi$ is a continuous function on the compact space $X,$ there exists a constant $C$ such that $$\psi\leq \varphi+C.$$ Thus by Lemma \ref{monotone2} and Lemma \ref{constantgo} we have that $$c[\psi]\leq c[\varphi+C]=c[\varphi]-C.$$ By symmetry we see that $c[\psi]-c[\varphi]$ is bounded on $\Delta(L)^{\circ}.$
\end{proof}

For Okounkov bodies we have that $$\Delta(mL)=m\Delta(L),$$ see e.g. \cite{Lazarsfeld}. The Chebyshev transforms also exhibit a homogeneity property.

\begin{proposition} \label{homo}
Let $\psi$ be a continuous metric on $L.$ Consider the metric $m\psi$ on $mL.$ For any $p\in \Delta(L)^{\circ}$ it holds that $$c[m\psi](mp)=mc[\psi](p).$$ 
\end{proposition}

\begin{proof}
We observe that trivially $A_{m\alpha,k}=A_{\alpha,km},$ as affine subspaces of $H^0(kmL),$ and hence $$F[m\psi](km\alpha,k)=F[\psi](km\alpha,km).$$ Let $\alpha(k)\to p\in \Delta(L)^{\circ}.$We get that
\begin{eqnarray*}
c[m\psi](mp)=|(mp,1)|c[m\psi](\frac{(mp,1)}{|(mp,1)|})=|(mp,1)|\lim_{k\to \infty}\frac{F[m\psi](km\alpha(k),k)}{k|(m\alpha(k),1)|}=\\=\lim_{k\to \infty}\frac{F[\psi](km\alpha(k),km)}{k}=mc[\psi](p).
\end{eqnarray*} 
\end{proof}

\section{The metric volume} \label{sectionmetric}

Recall that the volume of a line bundle $L$ is defined as $$\textrm{vol}(L):=\limsup_{k \to \infty}\frac{n!}{k^n}\textrm{dim}(H^0(kL)).$$ Let $(L,\psi)$ be a metrized line bundle. A metric version of the volume was introduced by Berman-Boucksom in \cite{Berman}. In \cite{Berman} it was called the \emph{energy at equilibrium}, but in this paper we simply call it the \emph{metric volume} because of the similarities it has with the ordinary volume function.

Given a metric $\psi$ one has a natural norm on the the spaces of holomorphic sections $H^0(kL),$ namely the supremum norm $$||s||_{k\psi,\infty}:=\sup\{|s(x)|e^{-k\psi(x)/2}: x\in X\}.$$ Let $\mathcal{B}^{\infty}(k\psi)\subseteq H^0(kL)$ be the unit ball with respect to this norm.

$H^0(kL)$ is a vector space, thus given a basis we can calculate the volume of $\mathcal{B}^{\infty}(k\psi)$ with respect to the associated Lebesgue measure. This will however depend on the choice of basis. But given a reference metric $\varphi$ one can compute the quotient  $$\frac{\textrm{vol}(\mathcal{B}^{\infty}(k\psi))}{\textrm{vol}(\mathcal{B}^{\infty}(k\varphi))}$$ and this will be independent of the choice of basis. The $k$:th $\mathcal{L}$-bifunctional is defined as $$\mathcal{L}_k(\psi,\varphi):=\frac{n!}{2k^{n+1}}\log \left(\frac{\textrm{vol}(\mathcal{B}^{\infty}(k\psi))}{\textrm{vol}(\mathcal{B}^{\infty}(k\varphi))}\right).$$

\begin{definition}
The metric volume of a metrized line bundle $(L,\psi),$ denoted by \\ $\textrm{vol}(L,\psi,\varphi),$ is defined as the limit 
\begin{equation} \label{introlimit2}
\textrm{vol}(L,\psi,\varphi):=\lim_{k\to \infty}\mathcal{L}_k(\psi,\varphi).
\end{equation}
\end{definition}

The metric volume obviously depends on the choice of $\varphi$ as a reference metric. But one can easily check that the difference of metric volumes $\textrm{vol}(L,\psi,\varphi)-\textrm{vol}(L,\psi',\varphi)$ is independent of the choice of reference (see \cite{Berman}), so without an explicit reference the metric volume is well-defined only up to a constant.

The definition of the metric volume is clearly reminiscent of the definition of the volume of a line bundles.

Consider the case where we let $\psi=\varphi+1.$ Since $$||s||^2_{k(\varphi+1)}=e^{-k}||s||^2_{k\varphi}$$ we get that $$\mathcal{B}^{\infty}(k(\varphi+1))=e^{-k}\mathcal{B}^{\infty}(k\varphi)$$ and thus by the homogeneity of the Lebesgue volume $$\mathcal{L}_k(\varphi+1,\varphi)=\frac{n!}{2k^{n+1}}\log e^{k2N_k},$$ where $$N_k=\dim_{\mathbb{C}}H^0(kL).$$ Thus the right hand side is equal to $$\frac{n!\dim_{\mathbb{C}}H^0(kL)}{k^n},$$ which converges to the volume of $L$ by definition.

In \cite{Berman} Berman-Boucksom prove that the limit (\ref{introlimit}) exists. They do this by proving that it actually converges to a certain integral over the space $X$ involving mixed Monge-Amp\`ere measures related to the metrics. This will be described in Section \ref{sectionenergy}. 

We now state our main result. 

\begin{theorem} \label{main}
Let $\psi$ and $\varphi$ be continuous metrics on $L.$ Then it holds that 
\begin{equation} \label{mainformula}
\textrm{vol}(L,\psi,\varphi)=n!\int_{\Delta(L)^{\circ}}(c[\varphi]-c[\psi])d\lambda,
\end{equation}
where $d\lambda$ denotes the Lebesgue measure on $\Delta(L)^{\circ}.$ 
\end{theorem}

The proof of Theorem \ref{main} will depend on the fact that one can also use $L^2$-norms to compute the Chebyshev transform of a continuous metric. This will be explained in the next section. It also yields a new proof of the existence of the limit (\ref{introlimit2}). 

\section{Bernstein-Markov norms}

\begin{definition}
Let $\mu$ be a positive measure on $X,$ and $\psi$ a continuous metric on a line bundle $L.$ One says that $\mu$ satisfies the Bernstein-Markov property with respect to $\psi$ if for each $\varepsilon > 0$ there exists $C=C(\varepsilon)$ such that for all non-negative $k$ and all holomorphic sections $s\in H^0(kL)$ we have that 
\begin{equation} \label{Bernstein}
\sup_{x\in X}\{|s(x)|^2e^{-k\psi(x)}\}\leq Ce^{\varepsilon k}\int_X|s|^2e^{-k\psi}d\mu.
\end{equation}
\end{definition}

If $\psi$ is a continuous metric on $L$ and $\mu$ a Bernstein-Markov measure on $X$ with respect to $\psi,$ we will call the $L^2$-norm on $H^0(kL)$ defined by $$||s||^2_{k\psi,\mu}:=\int_X|s|^2e^{-k\psi}d\mu$$ a Bernstein-Markov norm. We will also call the pair $(\psi,\mu)$ a Bernstein-Markov pair on $(X,L).$ 

For any continuous metric $\psi$ on $L$ there exist measures $\mu$ such that $(\psi,\mu)$ is a Bernstein-Markov pair. In fact any smooth volume form $dV$ on $X$ satisfies the Bernstein-Markov property with respect to any continuous metric, see Lemma 3.2 in \cite{Berman}.





We want to be able to use a Bernstein-Markov norm instead of the supremum norm to calculate the Chebyshev transform of a continuous metric $\psi.$

We pick a positive measure $\mu$ with the Bernstein-Markov property with respect to $\psi.$ For all $(k\alpha,k) \in \Gamma(L),$ let $t_{\alpha,k}$ be the section in $H^0(kL)$ of the form $$z^{k\alpha}+ \textrm{higher order terms}$$ that minimizes the $L^2$-norm $$||t_{\alpha,k}||^2_{k\psi,\mu}:=\int_{X}|t_{\alpha,k}|^2e^{-k\psi}d\mu.$$ It follows that $$<t_{\alpha,k},t_{\beta,k}>_{k\psi}=0$$ for $\alpha \neq \beta,$ since otherwise the sections $t_{\alpha,k}$ would not be minimizing. Hence $$\{t_{\alpha,k}: \alpha\in \Delta_k(L)\}$$ is an orthogonal basis for $H^0(kL)$ with respect to $||.||_{k\psi,\mu}.$ Indeed they are orthogonal, and by Lemma \ref{points} we have that $$\#\{t_{\alpha,k}: \alpha\in \Delta_k(L)\}=\#\Delta_k(L)=\textrm{dim}(H^0(kL)),$$ therefore it must be a basis. 

\begin{definition}
We define the discrete Chebyshev transform $F[\psi,\mu]$ of $(\psi,\mu)$ on $\Gamma$ by $$F[\psi,\mu](k\alpha,k):=\ln ||t_{\alpha,k}||_{k\psi,\mu}^2.$$ We also denote $\frac{1}{k}F[\psi,\mu](k\alpha,k)$ by $c_k[\psi,\mu](\alpha)$.
\end{definition} 

We will sometimes write $c_k[\psi]$ when we mean $c_k[\psi,\mu],$ considering $\mu$ as fixed.

\begin{proposition} \label{propro}
For any sequence $(k\alpha(k),k)$ in $\Gamma(L),$ $k\to \infty,$ such that $$\lim_{k\to \infty}\alpha(k)=p\in \Delta(L)^{\circ},$$ it holds that $$c[\psi](p)=\lim_{k\to \infty}c_k[\psi,\mu](\alpha(k)).$$
\end{proposition}

\begin{proof}
For a point $(k\alpha,k)\in \Gamma,$ let $t_{\alpha,k}$ be the minimizer with respect to the Bernstein-Markov norm. By the Bernstein-Markov property we get that $$||t_{\alpha,k}||^2_{\sup}\leq Ce^{\varepsilon k}||t_{\alpha,k}^{\mu}||^2_{\mu},$$ and hence 
\begin{equation} \label{hurra}
F[\psi](k\alpha,k)\leq F[\psi,\mu](k\alpha,k)+\ln C+\varepsilon k.
\end{equation}
Let $s$ be any section in $A_{\alpha,k}.$  We have that by definition $$||t_{\alpha,k}||^2_{\mu}\leq||s||^2_{\mu}\leq \mu(X) ||s||^2_{\sup},$$ so 
\begin{equation} \label{vabra}
F[\psi,\mu](k\alpha,k)\leq F[\psi](k\alpha,k)+\ln \mu(X).
\end{equation}
Equations (\ref{hurra}) and (\ref{vabra}) put together gives that 
\begin{equation} \label{eqproeq}
F[\psi](k\alpha,k)-\ln C-\varepsilon k\leq F[\psi,\mu](k\alpha,k)\leq F[\psi](k\alpha,k)+\ln \mu(X).
\end{equation}
It follows that $$\lim_{k\to \infty}\frac{F[\psi,\mu](k\alpha(k),k)}{k}=\lim_{k\to \infty}\frac{F[\psi](k\alpha(k),k)}{k}=c[\psi](p),$$ which gives the proposition.
\end{proof}

\begin{lemma} \label{constantgo2}
Let $\psi$ be a continuous metric on $L$ and consider the continuous metric on $L$ given by $\psi+C$ for some constant $C.$ Then it holds that $$F[\psi+C,\mu](k\alpha,k)=F[\psi,\mu](k\alpha,k)-kC.$$
\end{lemma}

\begin{proof}
This follows exactly as in the case of the suprumum norm, see proof of Lemma \ref{constantgo}.
\end{proof}

\begin{proposition} \label{monotone}
Let $(\psi,\mu)$ and $(\varphi,\nu)$ be two Bernstein-Markov pairs, and assume that $$\psi \leq \varphi$$ Then for every $\varepsilon>0$ there exists a constant $C'$ such that $$F[\psi,\mu](k\alpha,k)\geq F[\varphi,\nu](k\alpha,k)-C'-\varepsilon k.$$
\end{proposition}

\begin{proof}
Let $t_{\alpha,k}^{\psi}$ and $t_{\alpha,k}^{\varphi}$ be the minimizing sections with respect to the Bernstein-Markov norms $||.||_{k\psi,\mu}$ and $||.||_{k\varphi}$ respectively. From equation (\ref{eqproeq}) and Proposition \ref{monotone} we get that 
\begin{eqnarray*}
F[\psi,\mu](k\alpha,k)\geq F[\psi](k\alpha,k)-\ln C-\varepsilon k\geq  F[\varphi](k\alpha,k)-\ln C-\varepsilon k\geq \\ \geq F[\varphi,\nu]-\ln \nu(X)-\ln C-\varepsilon k.
\end{eqnarray*}
\end{proof}

\begin{proposition} \label{bounded}
For any two Bernstein-Markov pairs on $(X,L),$ $(\psi,\mu)$ and $(\varphi,\nu)$ the difference of the discrete Chebyshev transforms $$c_k[\psi,\mu]-c_k[\varphi,\nu]$$ is uniformly bounded on $\Delta(L)^{\circ}.$ 
\end{proposition}

\begin{proof}
By symmetry it suffices to find an upper bound. Let $\tilde{C}$ be a constant such that $\psi\leq \varphi+\tilde{C}.$ By Lemma \ref{constantgo2} and Proposition \ref{monotone} we get that 
\begin{eqnarray*}
c_k[\psi,\mu](\alpha)=\frac{1}{k}F[\psi,\mu](k\alpha,k)\geq \frac{1}{k}F[\varphi+C,\nu](k\alpha,k)-\frac{C'}{k}-\varepsilon=\\=\frac{1}{k}F[\varphi,\nu](k\alpha,k)-C-\frac{C'}{k}-\varepsilon=c_k[\varphi,\nu](\alpha)-C-\frac{C'}{k}-\varepsilon.
\end{eqnarray*}
The proposition follows.
\end{proof}

Finally let us consider $\mathcal{L}$ functionals using Bernstein-Markov norms instead of supremum norms.

Let $\mathcal{B}^2(\mu,k\varphi)$ denote the unit ball in $H^0(kL)$ with respect to the norm $\int_X |.|^2e^{-k\varphi}d\mu,$ i.e. $$\mathcal{B}^2(\mu,k\varphi):=\{s\in H^0(kL): \int_X |s|^2e^{-k\varphi}d\mu\leq 1\}.$$ Consider the quotient of the volume of two unit balls $$\frac{\textrm{vol} \mathcal{B}^2(\mu,k\varphi)}{\textrm{vol} \mathcal{B}^2(\nu,k\psi)}$$ with respect to the Lebesgue measure on $H^0(kL),$ where we by some linear isomorphism identify $H^0(kL)$ with $\mathbb{C}^N,$ $N=h^0(kL).$ In fact the quotient of the volumes does not depend on how we choose to represent $H^0(kL).$

\begin{lemma} \label{volumelemma}
\begin{equation} \label{volume}
\frac{\textrm{vol} \mathcal{B}^2(\mu,k\varphi)}{\textrm{vol} \mathcal{B}^2(\nu,k\psi)}=\frac{\textrm{det}(\int s_i \bar{s}_je^{-k\psi}d\nu)_{ij}}{\textrm{det}(\int s_i \bar{s}_je^{-k\varphi}d\mu)_{ij}},
\end{equation} 
where $\{s_i\}$ is any basis for $H^0(kL).$
\end{lemma}

\begin{proof}
First we show that the right hand side does not depend on the basis. Let $\{t_i\}$ be some orthonormal basis with respect to $\int |.|^2e^{-k\psi}d\nu$, and let $A=(a_{ij})$ be the matrix such that $$s_i=\sum a_{ij}t_j.$$ Then we see that 
\begin{equation} \label{change}
\int s_i \bar{s}_je^{-k\psi}d\nu=\int (\sum a_{ik}t_k)(\overline{\sum a_{jl}t_l})e^{-k\psi}d\nu=\sum a_{ik}\bar{a}_{jk}.
\end{equation}
Therefore by linear algebra we get that
\begin{equation}  \label{det}
\textrm{det}\left(\int s_i \bar{s}_je^{-k\psi}d\nu\right)_{ij}=\textrm{det}(AA^*)=|\textrm{det}A|^2.
\end{equation}
If we let $\{s_i'\}$ be a new basis, $$s_i'=\sum b_{ij}s_j, \qquad B=(b_{ij}),$$ then $$\textrm{det}\left(\int s_i' \bar{s}_j'e^{-k\psi}d\nu\right)_{ij}=|\textrm{det}B|^2\textrm{det}\left(\int s_i \bar{s}_je^{-k\psi}d\nu\right)_{ij}.$$ Since $|\textrm{det}B|^2$ also will show up in the denominator, we see that the quotient does not depend on the choice of basis.

Let as above $\{t_i\}$ be an orthonormal basis with respect to $\int |.|^2e^{-k\psi}d\nu$ and let $\{s_i\}$
be an orthonormal basis with respect to $\int |.|^2e^{-k\varphi}d\mu$ and let $$s_i=\sum a_{ij}t_j, \qquad{} A=(a_{ij}).$$ It is clear that $$\frac{\textrm{vol} \mathcal{B}^2(\mu,k\varphi)}{\textrm{vol} \mathcal{B}^2(\nu,k\psi)}=|\textrm{det}A|^2.$$ Note that the square in the right-hand side comes from the fact that we take the determinant of $A$ as a complex matrix. By equations (\ref{change}) and (\ref{det}) we also have that $$\textrm{det}\left(\int s_i \bar{s}_je^{-k\psi}d\nu\right)_{ij}=|\textrm{det}A|^2,$$ and since $\{s_i\}$ were chosen to be orthonormal $$\textrm{det}\left(\int s_i \bar{s}_je^{-k\varphi}d\mu\right)_{ij}=1.$$ The lemma follows.
\end{proof}

\begin{definition}
Let $(\varphi,\mu)$ and $(\psi,\nu)$ be two Bernstein-Markov pairs on $(X,L).$ The $L^2$-versions of the Donaldson $\mathcal{L}$ bifunctional, denoted by $\mathcal{L}_{k,2},$ is defined as $$\mathcal{L}_{k,2}(\varphi,\psi):=\frac{n!}{2k^{n+1}}\ln\left(\frac{\textrm{vol} \mathcal{B}^2(\mu,k\varphi)}{\textrm{vol} \mathcal{B}^2(\nu,k\psi)}\right).$$
\end{definition}

Let us to avoid confusion here denote the $k$:th $\mathcal{L}$ bifunctional using the supremum norm by $\mathcal{L}_{k,\infty}$

\begin{lemma} \label{twolimits}
For any two Bernstein-Markov pairs $(\varphi,\mu)$ and $(\psi,\nu)$ we have that $$\lim_{k\to \infty}\mathcal{L}_{k,2}(\varphi,\psi)=\lim_{k\to \infty}\mathcal{L}_{k,\infty}(\varphi,\psi)$$
if either limit exists.
\end{lemma}

\begin{proof}
By the Bernstein-Markov property we get that for any $\varepsilon>0$ there exists a constant $C$ such that $$C^{-1}e^{-k\varepsilon}\mathcal{B}^2(\mu,k\varphi)\subseteq \mathcal{B}^{\infty}(k\varphi)\subseteq ||\mu||\mathcal{B}^2(\mu,k\varphi).$$ Let $N_k$ denote the complex dimension of $H^0(kL).$ Because the Lebesgue volume on $H^0(kL)$ is $2N_k$-homogeneous we get that $$(C^{-1}e^{-k\varepsilon})^{2N_k}\textrm{vol}(\mathcal{B}^2(\mu,k\varphi))\leq \textrm{vol}(\mathcal{B}^{\infty}(k\varphi))\subseteq ||\mu||^{2N_k}\textrm{vol}(\mathcal{B}^2(\mu,k\varphi)).$$ Since $N_k\leq C'k^n$ for some constant $C'$ and in the expression for $\mathcal{L}_k$ we divide the logarithm of the volume by $k^{-n-1}$ we get that the two $\mathcal{L}$ functionals are asymptotically equal.
\end{proof}

\section{Proof of main theorem}

Here follows the proof of Theorem \ref{main}.

\begin{proof}
We let $\{s_i\}$ be a basis  for $H^0(kL)$ such that $$s_i=z^{k\alpha_i}+ \textrm{higher order terms},$$
where $\alpha_i \in \Delta_k(L)$ is some ordering of $\Delta_k(L).$ Let $$s_i=\sum a_{ij}t_{\alpha_j,k}^{\psi}, \qquad{} A=(a_{ij}).$$ From the proof of Lemma \ref{volumelemma} we see that
\begin{eqnarray*}
\textrm{det}\left(\int_X s_i\bar{s}_je^{-k\psi}d\nu\right)_{ij}=|\textrm{det}A|^2\textrm{det}\left(\int_X t_{\alpha_i,k}^{\psi}\bar{t}_{\alpha_j,k}^{\psi}e^{-k\psi}d\nu\right)_{ij}=\\=|\textrm{det}A|^2\prod_{\alpha\in \Delta_k(L)}||t_{\alpha,k}^{\psi}||^2,
\end{eqnarray*}
since $t_{\alpha,k}^{\psi}$ constitute an orthogonal basis. Also since the lowest term of $s_i$ is $z^{k\alpha_i}$ we must have that $a_{ij}=0$ for $j<i$ and $a_{ii}=1.$ Hence $\textrm{det}A=1,$ and consequently $$\textrm{det}\left(\int_X s_i\bar{s}_je^{-k\psi}d\nu\right)_{ij}=\prod_{\alpha\in \Delta_k(L)}||t_{\alpha,k}^{\psi}||^2.$$ From equation (\ref{volume}) we get that $$\mathcal{L}_{k,2}(\varphi,\psi)=\frac{n!}{k^n}\sum_{\alpha\in \Delta_k(L)}(c_k[\psi](\alpha)-c_k[\varphi](\alpha)).$$

For all $k$ let $\tilde{c}_k[\psi]$ denote the function on $\Delta(L)^{\circ}$ assuming the value of $c_k[\psi]$ in the nearest lattice point of $\Delta_k(L)$ (or the mean of the values if there are multiple lattice points at equal distance). Then $$\frac{n!}{k^n}\sum_{\alpha\in \Delta_k(L)}(c_k[\psi](\alpha)-c_k[\varphi](\alpha))=n!\int_{\Delta(L)^{\circ}}(\tilde{c}_k[\psi]-\tilde{c}_k[\varphi])d\lambda+\epsilon(k),$$ where the error term $\epsilon(k)$  goes to zero as $k$ tends to infinity since by Khovanskii's theorem we have that $\Delta_k(L)$ fills out more and more of $\Delta(L)^{\circ}\cap ((1/k)\mathbb{Z})^n$. By Propositions \ref{propro} and \ref{bounded} we can thus use dominated convergence to conclude that $$\lim_{k \to \infty}\mathcal{L}_{k,2}(\varphi,\psi)=n!\int_{\Delta(L)^{\circ}}(c[\psi]-c[\varphi])d\lambda.$$ Combined with Lemma \ref{twolimits} this proves the theorem. 
\end{proof}

\section{The Monge-Amp\`ere energy} \label{sectionenergy}

In \cite{Berman} Berman-Boucksom prove that the limit (\ref{introlimit}) exists. They do this by proving that it actually converges to a certain integral over the space $X$ involving mixed Monge-Amp\`ere measures related to the metrics. In order to describe this we need to introduce some concepts in pluripotential theory.

One can define a partial order on the space of metrics to a given line bundle. Let $\psi <_w \varphi$ if $$\psi \leq \varphi + O(1)$$ on $X.$ If a metric is maximal with respect to the order $<_w$, it is said to have minimal singularities. It is a fact that a metric with minimal singularities on a big line bundle is locally bounded on a dense Zariski-open subset of $X,$ see Section 1.4 in \cite{BEGZ}. On an ample line bundle, the metrics with minimal singularities are exactly those who are locally bounded.

Let $\psi$ and $\varphi$ be two locally bounded psh-metrics. By $\textrm{MA}_m(\psi,\varphi)$ we will denote the positive current $$\sum_{j=0}^m(dd^c\psi)^j\wedge (dd^c\varphi)^{m-j},$$ and by $\textrm{MA}(\psi)$ we will mean the positive measure $(dd^c \psi)^n.$ 

\begin{definition}
If $\psi$ and $\varphi$ are two psh metrics with minimal singularities, then we define the Monge-Amp\`ere energy of $\psi$ with respect to $\varphi$ as $$\mathcal{E}(\psi,\varphi):=\frac{1}{n+1}\int_\Omega (\psi-\varphi)\textrm{MA}_n(\psi,\varphi),$$ where $\Omega$ is a dense Zariski-open subset of $X$ on which $\psi$ and $\varphi$ are locally bounded. 
\end{definition}
 
\begin{remark}
In \cite{Berman} Berman-Boucksom use the notation $\mathcal{E}(\psi)-\mathcal{E}(\varphi)$ for what we denote by $\mathcal{E}(\psi,\varphi).$ Thus they consider $\mathcal{E}(\psi)$ as a functional defined only up to a constant.
\end{remark} 

An important aspect of the Monge-Amp\`ere energy (and a motivation for calling it an energy) is its cocycle property, i.e. that  $$\mathcal{E}(\psi,\varphi)+\mathcal{E}(\varphi,\psi')+\mathcal{E}(\psi',\psi)=0$$ for all metrics $\psi, \varphi$ and $\psi'.$ This is a reformulation of Corollary 4.2 in \cite{Berman}.

\begin{definition}
If $\psi$ is a continuous metric and $K$ a compact subset of $X$, the psh envelope of $\psi$ with respect to $K$, $P_K(\psi),$ is given by $$P_K(\psi):=\sup\{\varphi: \varphi \textrm{ psh metric on L}, \varphi\leq \psi \textrm{ on K}\}.$$
\end{definition}

For any $\psi$ and $K,$ as one may check, $P_K(\psi)$ will be psh and have minimal singularities. When $K=X,$ we will simply write $P(\psi)$ for $P_X(\psi).$

If $\psi$ and $\varphi$ are continuous metrics one can consider the composed functional $\mathcal{E}\circ P:$ $$\mathcal{E}\circ P(\psi,\varphi):=\frac{1}{n+1}\int_\Omega (P(\psi)-P(\varphi))\textrm{MA}_n(P(\psi),P(\varphi)).$$ 

We refer the reader to \cite{BEGZ} for a more thorough exposition on Monge-Amp\`ere measures and psh envelopes.

Theorem A in \cite{Berman} states that for Bernstein-Markov pairs the Donaldson $\mathcal{L}_k$ bifunctional converges to the composed functional $\mathcal{E}\circ P.$ Combined with our main result it yields the formula
\begin{equation} \label{superformula}
\frac{1}{n+1}\int_\Omega (P(\psi)-P(\varphi))\textrm{MA}_n(P(\psi),P(\varphi))=n!\int_{\Delta(L)^{\circ}}(c[\varphi]-c[\psi])d\lambda.
\end{equation}

If $\psi$ and $\varphi$ happen to be psh then we get that 
\begin{equation} \label{superformula2}
\mathcal{E}(\psi,\varphi)=n!\int_{\Delta(L)^{\circ}}(c[\varphi]-c[\psi])d\lambda.
\end{equation}

\section{Previous results}

Some instances of formula (\ref{mainformula}) are previously known. Here follows three such instances.  

\subsection{The volume as a metric volume}

In Section \ref{sectionmetric} we observed that for any metrized line bundle $(L,\psi)$ we have that 
\begin{equation} \label{previousequation1}
\textrm{vol}(L,\psi+1,\psi)=\textrm{vol}(L).
\end{equation}

Any minimizing section with respect to $\int |.|^2e^{-k\psi}$ will also minimize the norm $$\int |.|^2e^{-k(\psi+1)}=\int |.|^2e^{-k\varphi}.$$It follows that $c[\psi]-c[\varphi]$ is identically one. Therefore 
\begin{equation} \label{nofantasy}
\int_{\Delta(L)^{\circ}}(c[\psi]-c[\varphi])d\lambda=\textrm{vol}_{\mathbb{R}^n}(\Delta(L)).
\end{equation}
Equations (\ref{fantasy}) and (\ref{nofantasy}) and Theorem \ref{main} then gives us that $$\textrm{vol}(L)=n!\textrm{vol}_{\mathbb{R}^n}(\Delta(L)).$$ We have thus recovered Theorem A in \cite{Lazarsfeld}.

\subsection{Chebyshev constants and the transfinite diameter}

Let $K$ be a regular compact set in $\mathbb{C}.$ We let $||.||_K$ denote the norm which takes the supremum of the absolute value on $K.$ Let $P_k$ denote the space of polynomials in $z$ with $z^k$ as highest degree term. Let for any $k$ $$Y_k(K):=\inf\{||p||_K : p\in P_k\}.$$ One defines the Chebyshev constant $C(K)$ of $K$ as the following limit $$C(K):=\lim_{k\to \infty}(Y_k(K))^{1/k}.$$

Let $\{x_i\}_{i=1}^k$ be a set of $k$ points in $K.$ Let $d_k(\{x_i\})$ denote the product of their mutual distances, i.e. $$d_k(\{x_i\}):=\prod_{i<j}|x_i-x_j|.$$ One calls the points $\{x_i\}$ Fekete points if among the set of $k$-tuples of points in $K$ they maximize the function $d_k.$ Define $T_k(K)$ as $d_k(\{x_i\})$ for any set of Fekete points $\{x_i\}_{i=1}^k.$  Then the transfinite diameter $T(K)$ of $K$ is defined as $$T(K):=\lim_{k\to \infty}(T_k(K))^{1/\binom{k}{2}}.$$

We will now think of $\mathbb{C}$ as imbedded in the complex projective space $\mathbb{P}^1.$ Let $Z_0,Z_1$ be a basis for $H^0(\mathcal{O}(1)),$ therefore $[Z_0,Z_1]$ are homogeneous coordinates for $\mathbb{P}^1.$ Let $$z:=\frac{Z_1}{Z_0} \qquad{} \textrm{and} \qquad{} w:=\frac{Z_0}{Z_1}.$$ Let $p$ denote the point at infinity $$[0,1].$$ Then $w$ is a holomorphic coordinate around $p,$ and $Z_1$ is a local trivialization of the line bundle $\mathcal{O}(1)$ around $p.$ Thus we will identify a section $Z_0^{\alpha}Z_1^{k-\alpha}\in H^0(\mathcal{O}(k))$ with the polynomial $w^{\alpha}$ as well as with $z^{k-\alpha}.$ This means that the Okounkov body $\Delta(\mathcal{O}(1))$ of $\mathcal{O}(1)$ is the unit interval $[0,1]$ in $\mathbb{R}.$ We observe that a section $s\in H^0(\mathcal{O}(k))$ lies in $P_i$ as a polynomial in $z$ if and only if $$s=w^{k-i}+ \textrm{ higher order terms}.$$ For a section $s$ let $\tilde{s}$ denote the corresponding polynomial in $z.$ Consider the metric $P_K(\ln |Z_0|^2).$ It will be continuous since $K$ is assumed to be regular (see e.g. \cite{Berman}). Then we have the following lemma.

\begin{lemma} \label{lemmalemma}
For any $\alpha\in [0,1],$ i.e. that lies in the Okounkov body of $\mathcal{O}(1),$ we have that $$c[P_K(\ln |Z_0|^2)](\alpha)=2(1-\alpha)\ln C(K).$$
\end{lemma}

\begin{proof}
By basic properties of the projection operator $P_K$ (see \cite{Berman}) it holds that for for any section $s\in H^0(\mathcal{O}(k))$ 
\begin{equation} \label{rattita}
\sup_K\{|s|^2e^{-k\ln |Z_0|^2}\}=\sup_{\mathbb{P}^1}\{|s|^2e^{-kP_K(\ln |Z_0|^2)}\}.
\end{equation} 
Since the conversion to the $z$-variable means letting $Z_0$ be identically one, we also have that 
\begin{equation} \label{ohoh}
\sup_K\{|s|^2e^{-k\ln |Z_0|^2}\}=\sup_K\{|\tilde{s}|^2\}=||\tilde{s}||_K^2.
\end{equation}
We see that $s\in A_{\alpha,k}$ iff $\tilde{s}=z^{k-k\alpha}+\textrm{ lower order terms}.$ Hence $$F[P_K(\ln |Z_0|^2)](k\alpha,k)=2\ln Y_{k\alpha-k}(K),$$ and 
\begin{eqnarray*}
c[P_K(\ln |Z_0|^2)](\alpha)=\lim_{k\to \infty}\frac{F[P_K(\ln |Z_0|^2)](k\alpha,k)}{k}=\lim_{k\to \infty}\frac{2}{k}\ln Y_{k\alpha-k}(K)=\\=\lim_{k\to \infty}2(1-\alpha)\ln (Y_{k-k\alpha}(K))^{k-k\alpha}=2(1-\alpha)\ln C(K).
\end{eqnarray*}
\end{proof}

Let $K$ and $K'$ be two regular compact subsets of $\mathbb{C}.$ From Theorem \ref{main} and Lemma \ref{lemmalemma} we get that 
\begin{eqnarray*}
\mathcal{E}(P_{K'}(\ln |Z_0|^2),P_K(\ln |Z_0|^2))=\int_{(0,1)}(c[P_K(\ln |Z_0|^2)]-c[P_{K'}(\ln |Z_0|^2)])d\lambda(\alpha) \\ =\int_{(0,1)}\left(2(1-\alpha)\ln C(K)-2(1-\alpha)\ln C(K')\right)d\lambda(\alpha)=\ln C(K)-\ln C(K').
\end{eqnarray*}
On the other hand it follows from Corollary A in \cite{Berman} that 
\begin{equation} \label{brabra}
\ln T(K)-\ln T(K')=\mathcal{E}(P_{K'}(\ln |Z_0|^2),P_K(\ln |Z_0|^2)).
\end{equation}

Thus by Theorem \ref{main}, using Lemma \ref{lemmalemma} and equation (\ref{brabra}) we get that $$\ln T(K)-\ln T(K')=\ln C(K)-\ln C(K').$$ In fact it is easy to check that for the unit disc $D,$ $T(D)=C(D)=1,$ so we recover the classical result in potential theory that the transfinite diameter $T(K)$ and the Chebyshev constant $C(K)$ are equal.

For a thorough exposition on the subject of the transfinite diameter and capacities of compacts in $\mathbb{C}$ we refer the reader to the book \cite{Saff} by Saff-Totik.

\subsection{Invariant metrics on toric varieties}

Let $X$ be a smooth projective toric variety. We will view $X$ as a compactified $(\mathbb{C}^*)^n,$ such that the torus action on $X$ via this identification corresponds to the usual torus action on $(\mathbb{C}^*)^n.$ As is well-known, there is a polytope $\Delta$ naturally associated to the embedding $(\mathbb{C}^*)^n\subseteq X.$ We assume that $\Delta$ lies in the non-negative orthant of $\mathbb{R}^n.$ There is a line bundle $L_{\Delta}$ with a trivialization on $(\mathbb{C}^*)^n$ such that $$\Delta_k(L_{\Delta})=\Delta \cap (\frac{1}{k}\mathbb{Z})^n,$$ and any section $s\in H^0(kL_{\Delta})$ can in fact be written as a linear combination of the monomials $z^{\alpha}$ where $$\alpha\in k\Delta \cap \mathbb{Z}^n.$$

Let $dV$ be a smooth volume form on $X$ invariant under the torus action. Then it holds that for any torus invariant metric $\psi,$ $$\int_X z^{\alpha}\bar{z}^{\beta}e^{-k\psi}dV=0$$ when $\alpha \neq \beta.$ This follows from Fubini since trivially the monomials are orthogonal with respect to the Lebesgue measure on e.g. tori. Because of this for any torus invariant metric $\psi$ the minimizing sections $t_{a,k}^{\psi}$ are given by $z^{k\alpha},$ and consequently $$c_k[\psi,dV](\alpha)=\frac{1}{k}\ln \int_X |z^{k\alpha}|^{2}e^{-k\psi}dV.$$

Assume for simplicity that $\psi$ is positive.

\begin{lemma} \label{goodlemma}
For any strictly positive torus invariant metric $\psi$ we have that $$c[\psi](\alpha)=\ln \left(\sup_{z\in\mathbb{C}^n}\{|z^{\alpha}|^2e^{-\psi(z)}\}\right).$$
\end{lemma}    

\begin{proof}
We have that 
\begin{eqnarray*}
\int_X |z^{k\alpha}|^2e^{-k\psi}dV\leq dV(X)\sup_X \{|z^{k\alpha}|^2e^{-k\psi}\}=dV(X)\left(\sup_{z\in X} \{|z^{\alpha}|^2e^{-\psi(z)}\}\right)^k,
\end{eqnarray*}
which yieds the inequality $$c[\psi](\alpha)\leq \ln \left(\sup_{z\in X}\{|z^{\alpha}|^2e^{-\psi(z)}\}\right).$$ By the Bernstein-Markov property of $dV$ with respect to $\psi$ we get that 
\begin{eqnarray*}
\int_X |z^{k\alpha}|^2e^{-k\psi}dV\geq Ce^{-\varepsilon k}\sup_{z\in X} \{|z^{k\alpha}|^2e^{-k\psi(z)}\}=Ce^{-\varepsilon k}\left(\sup_{z\in X} \{|z^{\alpha}|^2e^{-\psi(z)}\}\right)^k.
\end{eqnarray*}
Using Proposition \ref{propro} it follows from this that $$c[\psi](\alpha)= \ln \left(\sup_{z\in X}\{|z^{\alpha}|^2e^{-\psi}\}\right).$$ Since $\psi$ is a metric on $L_{\Delta}$ it obeys certain growth conditions in $\mathbb{C}^n.$ In fact for $\alpha$ lying in the interior of $\Delta=\Delta(L_{\Delta})$ it holds that $$\sup_X\{|z^{\alpha}|^2e^{-\psi(z)}\}=\sup_{z\in \mathbb{C}^n}\{|z^{\alpha}|^2e^{-\psi(z)}\},$$ and the lemma follows.
\end{proof}

\begin{remark}
If we do not assume that the metric $\psi$ is strictly positive, the lemma still holds if we in the supremum replace $\psi$ with the projection $P(\psi)$.
\end{remark}

Let $\Theta$ denote the map from $(\mathbb{C}^*)^n$ to $\mathbb{R}^n$ that maps $z$ to $(\ln |z_1|,...,\ln |z_n|).$ Since we assumed $\psi$ to be torus invariant, the function $\psi \circ \Theta^{-1}$ is well-defined on $\mathbb{R}^n.$ We will denote $\psi \circ \Theta^{-1}$ by $\psi_{\Theta}.$ Since $\psi$ was assumed to be psh, it follows that $\psi_{\Theta}$ will be convex on $\mathbb{R}^n.$ Recall the definition of the Legendre transform. Given a convex function $g$ on $\mathbb{R}^n$ the Legendre transform of $g,$ denoted $g^*,$ evaluated in a point $p\in \mathbb{R}^n$ is given by $$g^*(p):=\sup_{x\in \mathbb{R}^n}\{\langle p,x \rangle -g(x) \}.$$ Observe that 
\begin{equation} \label{tatadidi}
\ln \left( (|z^{\alpha}|^2e^{-\psi})\circ \Theta^{-1}(x)\right)=2\langle \alpha,x \rangle-\psi_{\Theta}(x).
\end{equation}
Thus by equation (\ref{tatadidi}) and Lemma \ref{goodlemma} we get that $$c[\psi](\alpha)=2\left(\frac{\psi_{\Theta}}{2}\right)^*(\alpha).$$ Formula (\ref{superformula2}) now gives us that for any two invariant metrics $\psi$ and $\varphi$ on $L$ it holds that $$\mathcal{E}(\psi,\varphi)=2n!\int_{\Delta^{\circ}}\left(\frac{\varphi_{\Theta}}{2}\right)^*-\left(\frac{\psi_{\Theta}}{2}\right)^*d\lambda,$$ which is well-known in toric geometry. In fact this can be derived from the fact that the real Monge-Amp\`ere measure of a convex function is the pullback of the Lebegue measure with respect to the gradient of the convex function.   

\section{The Chebyshev transform on the zero-fiber}

Let us assume that $$z_1=0$$ is a local equation around $p$ for an irreducible variety which we denote by $Y.$  Let $H^0(X|Y,kL)$ denote the image of the restriction map from $H^0(X,kL)$ to $H^0(Y,kL_{|Y}),$  and let $\Gamma(X|Y,L)$ denote the semigroup $$\cup_{k\geq 0}\left(v(H^0(X|Y,kL))\times\{k\}\right) \subset \mathbb{N}^{n}.$$ Note that since $z_2,...,z_n$ are local coordinates on $Y,$ $v(H^0(X|Y,kL))$ will be a set of vectors in $\mathbb{N}^{n-1}.$

\begin{definition}
The restricted Okounkov body $\Delta_{X|Y}(L)$ is defined as the Okounkov body of the semigroup $\Gamma(X|Y,L).$
\end{definition}

\begin{lemma} \label{augmented}
If $Y$ is not contained in the augmented base locus $B_+(L),$ then $\Gamma(X|Y,L)$ generates $\mathbb{Z}^{n}$ as a group.
\end{lemma}

This is part of Lemma 2.16 in \cite{Lazarsfeld}.

\begin{remark}
The augmented base locus $B_+(L)$ of $L$ is defined as the base locus of any sufficiently small perturbation $L-\varepsilon A,$ where $A$ is some ample line bundle.
\end{remark}

Assume now that $Y$ is not contained in the augmented base locus $B_+(L).$ We will show that the Chebyshev transform $c[\psi]$ can be extended to the zero fiber, $$\Delta(L)_0:=\Delta(L)\cap \left(\{0\}\times \mathbb{R}^{n-1}\right),$$ in two different ways. 

From Theorem $4.24$ in \cite{Lazarsfeld} we get the following fact, 

\begin{equation} \label{restricted}
\Delta(L)_0=\Delta_{X|Y}(L).
\end{equation}

Note that since the Okounkov body lies in the positive orthant of $\mathbb{R}^n,$ $\Delta(L)_0$ is a part
of the boundary of $\Delta(L),$ hence the Chebyshev transform of a continuous metric is a priori not defined on the zero-fiber. Nevertheless, we want to show that one can extend the Chebyshev transform to the interior of zero-fiber $\Delta(L)_0$ in two different ways. 

First of all by restricting the discrete Chebyshev transform $F[\psi]$ on $\Gamma(L)$ to $\Gamma(X|Y,L)$ we get in the ordinary way a convex function on the interior of $\Delta(L)_0,$ thanks to  Lemma \ref{augmented} and (\ref{restricted}). We will call this function the restricted Chebyshev transform and denote it by $c_{X|Y}[\psi].$

It is not clear that $c_{X|Y}[\psi]$ gives a continuous extension of $c[\psi].$ However we will show that there is a continuous extension of $c[\psi]$ to the interior of the zero fiber which is at least bounded from above by $c_{X|Y}[\psi].$ 

To do this, we need to know how $\Gamma$ behaves near this boundary, something which Theorem \ref{khov} does not tell us anything about.

\begin{lemma} \label{fiberbound}
Assume Y is not contained in the augmented base locus of $L$, and let $p$ be any point in the interior of $\Delta(L)_0.$ Let $\Sigma^{\mathbb{Z}}_{n+1}$ denote the unit simplex in $\mathbb{Z}^{n+1},$ $\Sigma^{\mathbb{R}}_{n-1}$ the unit simplex in $\mathbb{R}^{n-1},$ and let $S$ denote the simplex  $\{0\}\times \Sigma^{\mathbb{R}}_{n-1}\times \{0\}.$ Then $\Gamma(L)$ contains a translated unit simplex $(\alpha,k)+\Sigma_{n+1}$ such that $(kp,k)$ lies in the interior of the $(n-1)$-simplex $$(\alpha,k)+S$$ (i.e interior with respect to the $\mathbb{R}^{n-1}$ topology). 
\end{lemma}

\begin{proof}
By Lemma \ref{augmented} we may use Lemma \ref{goodygoody} in combination with equation (\ref{restricted}) to reach the conclusion that for large $k,$ there are sections $s_k$ such that $(p,k)$ lies in the interior of  $(v(s_k),k)+S$ with respect to the $\mathbb{R}^{n-1}$ topology. We may write $L$ as a difference of two very ample divisors $A$ and $B.$ We may choose $B$ such that $\Delta_1(B)$ contains $\Sigma_n$ in $\mathbb{Z}^{n}$, and A such that $\Delta_1(A)$ contains origo. Now $$kL=B+(kL-B).$$ Since $L$ is big, for $k$ large we can find sections $s_k'\in H^0(kL-B)$ such that $v(s_k')=v(s_k).$  We get that $$(v(s_k),k)+\Sigma_n\subseteq \Gamma(L),$$ by multiplying $s_k'$ by the sections of $B$ corresponding to the points in the unit simplex $\Sigma_n\subseteq \Delta_1(B).$ Also observe that $$(k+1)L=A+(kL-B).$$ Now by multiplying $s_k'$ with the section of $A$ corresponding to origo in $\Delta_1(A)$ we get $$(v(s_k'),k)+(0,...,0,1)\subseteq \Gamma(L).$$ Since $$\Sigma_n\times \{0\}\cup (0,...,0,1)=\Sigma_{n+1}$$ we get $$(v(s_k'),k)+\Sigma_{n+1}\subseteq \Gamma(L).$$
\end{proof}

\begin{remark}
The proof is very close to the argument in \cite{Lazarsfeld} which shows the existence of a unit simplex in $\Gamma(L),$ when $L$ is big. The difference here is that we need to control the position of the unit simplex, but the main trick of writing $L$ as a difference of two very ample divisors is the same.
\end{remark}

\begin{lemma} \label{khov23}
Let $p$ be as in the statement of Lemma \ref{fiberbound}. Then there exists a neighbourhood $U$ of $p$ such that if we denote the intersection $U\cap \Delta(L)$ by $\tilde{U},$ for $k$ large it holds that $$(k\tilde{U},k)\cap \mathbb{Z}^{n+1}\subseteq \Gamma(L).$$
\end{lemma} 

\begin{proof}
Let $(\alpha,m)+\Sigma^{\mathbb{Z}}_{n+1}\subseteq \Gamma(L)$ be as in the statement of Lemma \ref{fiberbound}, and let $D^{\mathbb{Z}}\subseteq \Gamma(L)$ denote the set $$D^{\mathbb{Z}}:=(\alpha,m)+\Sigma_n^{\mathbb{Z}}\times\{0\}=(\alpha+\Sigma_n^{\mathbb{Z}})\times\{m\}.$$ Let also $D^{\mathbb{R}}$ denote the set $$D^{\mathbb{R}}:=(\alpha+\Sigma_n^{\mathbb{R}})\times\{m\}.$$Since trivially $$\underbrace{\Sigma_n^{\mathbb{Z}}+...+\Sigma_n^{\mathbb{Z}}}_{k}=(k\Sigma_n^{\mathbb{R}})\cap \mathbb{Z}^n,$$ we have that $$(kD^{\mathbb{R}},km)\cap \mathbb{Z}^{n+1}=\underbrace{D^{\mathbb{Z}}+...+D^{\mathbb{Z}}}_{k} \subseteq \Gamma(L).$$ Therefore the lemma holds when $k$ is a multiple of $m.$ Furthermore, since $m$ and $m+1$ are relatively prime, if $k$ is greater than $m(m+1)$ we can write $$k=k_1m+k_2(m+1),$$ where both $k_1$ and $k_2$ are non-negative, and $k_2\leq m.$ Thus we consider the set $$\underbrace{D^{\mathbb{Z}}+...+D^{\mathbb{Z}}}_{k_1}+k_2(\alpha,m+1)\subseteq \Gamma(L).$$ Because of the bound $k_2\leq m,$ and since $(\alpha,m+1)$ lies on the zero fiber, for a neighbourhood $\tilde{U}$ of $p,$ when $k$ gets large we must have that $$(k\tilde{U},k)\cap \mathbb{Z}^{n+1}\subseteq \underbrace{D^{\mathbb{Z}}+...+D^{\mathbb{Z}}}_{k_1}+k_2(\alpha,m+1) \subseteq \Gamma(L).$$
\end{proof}

\begin{corollary} \label{corzero}
Assume $Y$ is not contained in the augmented base locus  of $L.$ Then the chebyshev function $c[\psi]$ has a continuous extension to the interior of the zero-fiber, $\Delta(L)_0,$ and it is continuous and convex on its extended domain $\Delta(L)^{\circ}\cup\Delta(L)_0^{\circ}.$
\end{corollary}

\begin{proof}
Using Lemma \ref{khov23} and the subadditivity of $F[\psi]$ yields that $c[\psi]$ is bounded in a neighbourhood of any point $p$ in the interior of $\Delta(L)_0$. It is an elementary fact that any convex function defined on an open half space which is locally bounded near the boundary has a convex continuous extension to the boundary. Therefore it follows that $c[\psi]$ has a convex continuous extension to the interior of $\Delta(L)_0.$
\end{proof}

\begin{lemma} \label{amplerestr}
Assume $L$ is ample, and $\psi$ is a continuous metric. Then for any regular compact set $K$ it holds that the projection $P_K(\psi)$ also is continuous. In particular, since $X$ is regular, $P(\psi)$ is continuous when $L$ is ample. 
\end{lemma}

\begin{proof}
See e.g. \cite{Berman}.
\end{proof}

We will have use for the Ohsawa-Takegoshi extension theorem. We choose to record one version (see e.g. \cite{Demailly}).

\begin{theorem} \label{ohsawa}
Let $L$ be a holomorphic line bundle and let $S$ be a divisor. Assume that $L$ and $S$ have metrics $\Psi_L$ and $\Psi_S$ respectively satisfying $$dd^c\Psi_L\geq(1+\delta)dd^c\Psi_S+dd^c\Psi_{K_X},$$ where $\Psi_{K_X}$ is some smooth metric on the canonical bundle $K_X.$ Assume also that $$dd^c\Psi_L \geq dd^c(\Psi_S+\Psi_{K_X}).$$ Then any holomorphic section $\tilde{t}$ of the restriction of $L$ to $S$ extends holomorphically to a section $t$ of $L$ over $X$ satisfying $$\int_X |t|^2 e^{-\Psi_L}\omega_n \leq C_{\delta}\int_S |\tilde{t}|^2e^{-\Psi_L}\frac{dS}{|ds|^2e^{-\Psi_S}}.$$ Here $\omega_n$ is a smooth volume form on $X$ and $dS$ is a smooth volume form on $S.$
\end{theorem}

\begin{lemma} \label{newnewlemma}
Suppose $L$ is ample. Let $A$ be an ample line bundle, with a holomorphic section $s$ such that locally $s=z_1$. Also assume that the zero-set of $s$, which we will denote by $Y$, is a smooth submanifold. Then for all $\alpha \in \Delta_{X|Y}(L)$ we have that
\begin{equation} \label{restrcheb}
c_{X|Y}[\varphi](\alpha)=c_Y[P(\varphi)_{|Y}](\alpha).
\end{equation}
\end{lemma}

\begin{proof}
We may choose $\tilde{z}_1=z_2,...,\tilde{z}_{n-1}=z_n$ as holomorphic coordinates on $Y$ around $p.$ We consider the discrete Chebyshev transforms of the restrictions of $P(\varphi)$ and $P(\psi)$ to $Y.$ Since $L$ is ample, by Lemma \ref{amplerestr} $P(\varphi)$ and $P(\psi)$ are continuous, therefore the restrictions will also be continuous psh-metrics on $L_{|Y},$ therefore the Chebyshev transforms $c_Y[P(\varphi)_{|Y}]$ and  $c_Y[P(\psi)_{|Y}]$ are well-defined.

We note that if $t\in H^0(X,kL)$ and $$t=z^{k(0,\alpha)}+\textrm{higher order terms},$$ the restriction of $t$ to $Y$ will be given by $$t_{|Y}=\tilde{z}^{k\alpha}+\textrm{higher order terms}.$$ Furthermore  
\begin{eqnarray*}
\sup_Y \{|t_{|Y}|^2e^{-kP(\varphi)}\}\leq \sup_X \{|t|^2 e^{-kP(\varphi)}\}.
\end{eqnarray*}

This gives the inequality $$c_{X|Y}[\varphi](\alpha) \geq c_Y[P(\varphi)_{|Y}](\alpha),$$ by taking $t$ to be some minimizing section with respect to the supremum norm on $X.$

For the opposite inequality we use Proposition \ref{propro} which says that one can use Bernstein-Markov norms to compute the Chebyshev transform.

If $\tilde{t}\in H^0(Y,kL_{|Y}),$ $$\tilde{t}=\tilde{z}^{k\alpha}+\textrm{higher order terms},$$ then if $k$ is large enough there exists a section $t\in H^0(X,kL)$ such that $t_{|Y}=\tilde{t}.$ This is because we assumed $L$ to be ample, so we have extension properties (by e.g. Ohsawa-Takegoshi). We observe that any such extension must look like $$t=z^{k(0,\alpha)}+\textrm{higher order terms},$$ because if we had that $$t=z^{k(\beta_1,\beta)}+\textrm{higher order terms}$$ with $\beta_1>0,$ then since all higher order terms also restrict to zero, $$t_{Y}=0,$$ which is a contradiction. 

Let $\Psi$ be some smooth strictly positive metric on $L.$ Then for some $m$ $$dd^cm\Psi>(1+\delta)dd^c\Psi_A+dd^c\Psi_{K_X}$$ and $$dd^cm\Psi>dd^c\Psi_A+dd^c\Psi_{K_X},$$ where $\Psi_A$ and $\Psi_{K_X}$ are metrics on $A$ and $K_X$ respectively. We have that $dd^cP(\varphi)\geq 0,$ hence $$dd^c((k-m)P(\varphi)+m\Psi)>(1+\delta)dd^c\Psi_A+dd^c\Psi_{K_X}$$ and $$dd^c((k-m)P(\varphi)+m\Psi)>dd^c\Psi_A+dd^c\Psi_{K_X}$$for all $k>m.$ Since $P(\varphi)$ is continuous hence locally bounded, we also have that for some constant $C,$ $$\Psi-C<P(\varphi)<\Psi+C.$$ We can apply Theorem \ref{ohsawa} to these metrics, and get that for large $k,$ given a $\tilde{t}\in H^0(Y,kL_{|Y})$ there exists an extension $t\in H^0(X,kL)$ such that 
\begin{eqnarray*}
\int_X |t|^2 e^{-kP(\varphi)}\omega_n \leq e^{mC}\int_X |t|^2 e^{-(k-m)P(\varphi)-m\Psi}d\mu \\ \leq e^{mC} C_{\delta} \int_Y |\tilde{t}|^2 e^{-(k-m)P(\varphi)-m\Psi}d\nu \leq e^{2mC}C_{\delta} \int_Y |\tilde{t}|^2 e^{-kP(\varphi)}d\nu,
\end{eqnarray*}
where $C_{\delta}$ is constant only depending on $\delta$ and $d\nu$ is a smooth volume form on $Y$.  
By letting $\tilde{t}$ be the minimizing section with respect to $\int_Y |.|^2e^{-kP(\varphi)}d\nu$ and using  Proposition \ref{propro} we get that $$c_{X|Y}[\varphi](\alpha) \leq c_Y[P(\varphi)_{|Y}](\alpha),$$ since $$\int_X |t|^2 e^{-k\varphi}\omega_n \leq \int_X |t|^2 e^{-kP(\varphi)}\omega_n.$$ 
\end{proof}

\begin{proposition} \label{restrenergy}
Let $L,$ $A$ and $Y$ be as in the statement of Lemma \ref{newnewlemma}. Then we have that $$\textrm{vol}(L_{|Y},P(\varphi)_{|Y},P(\psi)_{|Y})=(n-1)!\int_{\Delta(L)_0} (c_{X|Y}[\psi]-c_{X|Y}[\varphi])(\alpha)d\alpha.$$
\end{proposition}

\begin{proof}
The proposition follows from Lemma \ref{newnewlemma} by integration of equality (\ref{restrcheb}) over the interior of the zero-fiber, and Theorem \ref{main} which says that 
$$\textrm{vol}(L_{|Y},P(\varphi)_{|Y},P(\psi)_{|Y})=(n-1)!\int_{\Delta(L_{|Y})} c_Y[P(\psi)_{|Y}]-c_Y[P(\varphi)_{|Y}]d\lambda.$$
\end{proof}

We will cite Proposition 4.7 from \cite{Berman} which is a recursion formula relating the metric volume with the restricted version.

\begin{proposition} \label{recursion}
Suppose $L$ is ample, let $s\in H^0(L),$ and let $Y$ be the smooth submanifold defined by $s$. Let $\psi$ and $\varphi$ be two continuous metrics. Then
\begin{eqnarray*}
(n+1)\textrm{vol}(L,\psi,\varphi)-n\textrm{vol}(L_{|Y},P(\varphi)_{|Y},P(\psi)_{|Y})=\\=\int_X (\ln |s|^2-P(\varphi))\textrm{MA}(P(\varphi))-\int_X(\ln |s|^2-P(\psi))\textrm{MA}(P(\psi)).
\end{eqnarray*}
\end{proposition}

In particular, combining Theorem \ref{main}, Proposition \ref{restrenergy} and Proposition \ref{recursion} we get the following.

\begin{proposition} \label{kruit}
Let $L,$ $s$ and $Y$ be as in Proposition \ref{recursion}. Then it holds that 
\begin{eqnarray*}
\int_{\Delta(L)^{\circ}}(c_X[\varphi]-c_X[\psi])d\lambda_n=\frac{1}{n+1}\int_{\Delta(L)_0^{\circ}}(c_{X|Y}[\varphi]-c_{X|Y}[\psi])d\lambda_{n-1}+\\+\frac{1}{(n+1)!}\int_X (\ln |s|^2-P(\varphi))\textrm{MA}(P(\varphi))-\frac{1}{(n+1)!}\int_X(\ln |s|^2-P(\psi))\textrm{MA}(P(\psi)).
\end{eqnarray*}
\end{proposition}

\section{Directional Chebyshev constants in $\mathbb{C}^n$}

In \cite{Bloom} Bloom-Levenberg define the metrized version of the directional Chebyshev constants originally introduced by Zaharjuta in \cite{Zaharjuta}. In this section we will describe how this relates to the Chebyshev transforms we have introduced. 

The setting in \cite{Bloom} is as follows. Let $<_1$ be the order on $\mathbb{N}^n$ such that $\alpha <_1\beta$ if $|\alpha|<|\beta|,$ or if $|\alpha|= |\beta|$ and $\alpha <_{\textrm{lex}} \beta.$ Let $P_{\alpha}$ denote the set of polynomials $p(z_1,...,z_n)$ in the variables $z_i$ such that $$p=z^{\alpha}+\textrm{ lower order terms}.$$ Observe that here we want lower order terms, and not higher order terms. Let $K$ be a compact set and $h$ an admissible metric function on $K.$ For any $\alpha \in \mathbb{N}^n$ they define the metrized Chebyshev constant $Y_3(\alpha)$ as $$Y_3(\alpha):= \inf \{\sup_{z\in K} \{|h(z)^{|\alpha|}p(z)|\}: p\in P_{\alpha}\}.$$ They then show that the limit $$\tau^{h}(K,\theta):=\lim_{\alpha/\textrm{deg}(\alpha)\to \theta}Y_3(\alpha)^{1/\textrm{deg}(\alpha)}$$ exists. These limits are called directional Chebyshev constants. 

In our setting we wish to view $\mathbb{C}^n$ as an affine space lying in $\mathbb{P}^n.$ Also, polynomials in $z_i$ can be interpreted as sections of multiples of the line bundle $\mathcal{O}(1)$ on $\mathbb{P}^n$ in the following sense. 
Let $Z_0,...,Z_n$ be a basis for $H^0(\mathcal{O}(1))$ on $\mathbb{P}^n$, and identify them with the homogeneous coordinates $[Z_0,...,Z_n].$ We can choose $$p:=[1:0:...:0]$$  to be our base point, and let $z_i:=\frac{Z_i}{Z_0}$ be holomorphic coordinates around $p.$ We also let $Z_0$ be our local trivialization of the bundle. Given a section $s\in H^0(\mathcal{O}(k))$we  represent it as a function in $z_i$ by dividing by a power of $Z_0$ $$\frac{s}{Z_0^k}=\sum a_{\alpha}z^{\alpha}.$$ Therefore we see that $$Z^{(\alpha_0,\alpha_1,...,\alpha_n)} \mapsto z^{(\alpha_1,...,\alpha_n)}.$$ 

We could also choose a different set of coordinates. Let $$q:=[0:...:0:1]$$ be our new base point, and let $w_i:=\frac{Z_i}{Z_n}$ be coordinates around $q.$ Let $Z_n$ be the local trivialization around $q.$  
Given a section $s\in H^0(\mathcal{O}(k))$ we represent it as a function in $w_i$ by dividing by a power of $Z_n$ $$\frac{s}{Z_n^k}=\sum b_{\alpha}w^{\alpha}.$$ Hence $$Z^{(\alpha_0,\alpha_1,...,\alpha_n)} \mapsto w^{(\alpha_0,...,\alpha_{n-1})}.$$ To define Chebyshev transforms we need an additive order on $\mathbb{N}^n$. Since the semigroup $\Gamma(\mathcal{O}(1))$ will not depend on the order, we are free to choose any additive order. Let $<_2$ be the order which corresponds to inverting the order $<_1$ with respect to the $z_i$ variables, i.e. $$(\alpha_0,...,\alpha_{n-1})<_2(\beta_0,...,\beta_{n-1})$$ iff $$(\beta_1,...,\beta_n)<_1(\alpha_1,...,\alpha_n).$$ Therefore 
\begin{equation} \label{zw}
z^{(\alpha_1,...,\alpha_n)}+\textrm{ lower order terms}=w^{(\alpha_0,...,\alpha_{n-1})}+\textrm{ higher order terms}.
\end{equation} 

We may identify the metric function $h$ with a metric $h=e^{-\psi/2}$ on $\mathcal{O}(1).$ Consider the metric $P_K(\psi).$ For simplicity assume that $K$ is regular. Since $\mathcal{O}(1)$ is ample from Lemma \ref{amplerestr} it follows that $P_K(\psi)$ is continuous, therefore the Chebyshev transform $c[P_K(\psi)]$ is well-defined. It is a simple fact that
\begin{equation} \label{simple}
\sup_{z\in K}\{|s(z)|^2e^{-k\psi(z)}\}=\sup_{z\in\mathbb{P}^n}\{|s(z)|^2e^{-kP_K(\psi)(z)}\}.
\end{equation}

Let $\alpha_0=0,$ and let $k=\sum_1^n \alpha_i.$ By (\ref{zw}) we see that $s\in A_{(\alpha_0,...,\alpha_{n-1}),k}$ iff it is on the form $$z^{(\alpha_1,...,\alpha_n)}+\textrm{ lower order terms}.$$ By (\ref{simple}) it follows that $$\ln Y_3(\alpha_1,...,\alpha_n)=F[P_K(\psi)](k\alpha,k).$$ 
Thus we get that for $\theta=(\theta_1,...,\theta_n)\in \Sigma^0$
\begin{equation} \label{directed}
c_{\mathbb{P}^n|\mathbb{P}^{n-1}}[P_K(\psi)](\theta_1,...,\theta_{n-1})=2\ln \tau^h(\theta_1,...,\theta_n).
\end{equation}

Observe that the order $<_2$ we used to defined the Chebyshev transform has the property that $(0,\alpha)<_2(\beta_1,\beta)$ when $\beta_1>0.$ It was this property of the lexicographic order we used in the proof of Proposition \ref{restrenergy}. Therefore the theorem holds also for Chebyshev transforms defined using $<_2$ instead of $<_{\textrm{lex}}.$ Let $(K',h')$ be another metrized set in $\mathbb{C}^n$, and let $\psi'$ be the corresponding metric on $\mathcal{O}(1)$ associated to $h'.$ Then integrating (\ref{directed}) gives us that 
\begin{eqnarray} \label{trutru}
\frac{1}{\textrm{meas}(\Sigma^0)}\int_{\Sigma^0}\ln \tau^h(K,\theta)-\ln \tau^{h'}(K',\theta)d\theta= \nonumber \\=\frac{(n-1)!}{2}\int_{\Delta(\mathcal{O}(1))_0}c_{\mathbb{P}^n|\mathbb{P}^{n-1}}[P_K(\psi)]-c_{\mathbb{P}^n|\mathbb{P}^{n-1}}[P_{K'}(\psi')]d\theta,
\end{eqnarray}
where $Y:=\{Z_0=0\}.$ Here we used that $\Delta(\mathcal{O}(1))_0$ is a $(n-1)$-dimensional unit simplex, and thus $$\textrm{meas}(\Delta(\mathcal{O}(1))_0)=\frac{1}{(n-1)!}.$$  

Bloom-Levenberg define a metrized transfinite diameter $d^{h}(K)$ of $K$ which is given by $$d^{h}(K):=\textrm{exp}\left(\frac{1}{\textrm{meas}(\Sigma^0)}\int_{\Sigma^0}\ln \tau^h(K,\theta)d\theta\right).$$ There is also another transfinite diameter, $\delta^h(K)$, which is defined as a limit of certain Vandermonde determinants. By Corollary A in \cite{Berman} we have that $$\ln \delta^h(K)-\ln \delta^{h'}(K')=\frac{(n+1)}{2n}\mathcal{E}(P_{K'}(\psi'),P_K(\psi)).$$ Then by Theorem \ref{main}, equation (\ref{trutru}) and Proposition \ref{kruit} we get that
\begin{eqnarray*}
\ln \delta^h(K)-\ln \delta^{h'}(K')=\\=\ln d^{h}(K)-\ln d^{h'}(K')+\frac{1}{n}\int_{\mathbb{P}^n}\frac{1}{2}(\ln |Z_0|^2-P_K(\psi))\textrm{MA}(P_K(\psi))-\\-\frac{1}{n}\int_{\mathbb{P}^n}\frac{1}{2}(\ln |Z_0|^2-P_{K'}(\psi'))\textrm{MA}(P_{K'}(\psi')).
\end{eqnarray*}
In fact, the positive measure $\textrm{MA}(P_K(\psi))$ has support on $K,$ and $P_K(\psi)=\psi$ a.e. with respect to $\textrm{MA}(P_K(\psi))$. In the notation of \cite{Bloom}, $(\psi-\ln|Z_0|^2)/2$ is denoted $Q,$ and $\textrm{MA}(P_K(\psi))$ is denoted $(dd^cV_{K,Q}^*)^n.$ Thus in their notation
\begin{eqnarray*}
\ln \delta^h(K)-\ln \delta^{h'}(K')=\\=\ln d^{h}(K)-\ln d^{h'}(K')-\frac{1}{n}\int_K Q (dd^cV_{K,Q}^*)^n+\frac{1}{n}\int_{K'}Q'(dd^cV_{K',Q'}^*)^n.
\end{eqnarray*} 
For the unit ball $B,$ with $h\equiv 1 \equiv |Z_0|^2$ and therefore $Q_h=0,$ it is straight-forward to show that we have an equality $$\delta^h(B)=d^{h}(K).$$ Using this we get that $$\ln \delta^h(K)=\ln d^{h}(K)-\frac{1}{n}\int_K Q (dd^cV_{K,Q}^*)^n.$$ By taking the exponential we have derived the formula of Theorem 2.9 in \cite{Bloom}. 

\section{Chebyshev transforms of metrized $\mathbb{Q}$- and $\mathbb{R}$-divisors}

Because of the homogeneity of Okounkov bodies, one may define the Okounkov body $\Delta(D)$ of any big $\mathbb{Q}$-divisor $D.$ Set $$\Delta(D):=\frac{1}{p}\Delta(pD)$$ for any $p$ that clears the denominators in $D.$ In \cite{Lazarsfeld} Lazarsfeld-Musta\c{t}\u{a} show that this mapping of a $\mathbb{Q}$-divisor to its Okounkov body has a continuous extension to the class of $\mathbb{R}$-divisors.

In Proposition \ref{homo} we saw that Chebyshev transforms also are homogeneous under scaling. Therefore we may define the Chebyshev transform of a $\mathbb{Q}$-divisor $D$ with metric $\psi,$ by letting
\begin{equation} \label{qdivisor}
c[\psi](\alpha)=\frac{1}{p}c[p\psi](p\alpha), \qquad \alpha \in \Delta(D)^{\circ},
\end{equation}
for any $p$ clearing the denominators in $D.$ We wish to show that this can be extended continuously to the class of metrized $\mathbb{R}$-divisors.

We will use the construction introduced in \cite{Lazarsfeld}. Let $D_1,...,D_r$ be divisors such that every divisor is numerically equivalent to a unique sum $$\sum a_iD_i, \qquad{} a_i \in \mathbb{Z}.$$ Lazarsfeld-Musta\c{t}\u{a} show that for effective divisors the coefficients $a_i$ may be chosen non-negative.

\begin{definition}
The semigroup of $X,$ $\Gamma(X),$ is defined as $$\Gamma(X):=\bigcup_{a\in \mathbb{N}^r}\left(v(H^0(\mathcal{O}_X(\sum a_iD_i)))\times \{a\}\right)\subseteq \mathbb{Z}^{n+r},$$ where $v$ stands for the usual valuation, $$s=z^{\alpha}+ \textrm{ higher order terms} \qquad{} \Rightarrow \qquad{} s \mapsto \alpha.$$ 
\end{definition}

Lazarsfeld-Musta\c{t}\u{a} show in \cite{Lazarsfeld} that $\Gamma(X)$ generates $\mathbb{Z}^{n+r}$ as a group. 

Let $\Sigma(\Gamma(X))$ denote the closed convex cone spanned by $\Gamma(X),$ and let for $a \in \mathbb{N}^r$ $$\Delta(a):=\Sigma(\Gamma(X))\cap (\mathbb{R}^n\times \{a\}).$$

By \cite{Lazarsfeld} for any big $\mathbb{Q}$-divisor $D=\sum a_iD_i,$ $$\Delta(a)=\Delta(D), \qquad{} a=(a_1,...,a_r).$$

Let for each $1\leq i\leq r$ $\psi_i$ be a continuous metrics on $D_i.$ Then for $a\in \mathbb{N}^r,$ $\sum a_i \psi$ is a continuous metric on $\sum a_iD_i.$ For an element $(\alpha,a)\in \Gamma(X),$ let $A_{\alpha,a}\subseteq H^0(\sum a_iD_i)$ be the set of sections of the form $$z^{\alpha}+ \textrm{ higher order terms}.$$ 

\begin{definition}
The discrete global Chebyshev transform $F[\psi_1,...,\psi_r]$ is defined by $$F[\psi_1,...,\psi_r](\alpha,a):=\inf\{\ln ||s||^2_{\alpha,a}: s\in A_{\alpha,a}\}$$ for $(\alpha,a)\in \Gamma(X).$
\end{definition}

\begin{lemma} \label{soon}
$F[\psi_1,...,\psi_r]$ is subadditive on $\Gamma(X)$.
\end{lemma}

\begin{proof}
If $s\in H^0(\mathcal{O}_X(\sum a_iD_i)),$ $$s=z^{\alpha}+ \textrm{ higher order terms},$$ and $t\in H^0(\mathcal{O}_X(\sum b_iD_i)),$ $$t=z^{\beta}+ \textrm{ higher order terms},$$ then $st\in H^0(\mathcal{O}_X(\sum (a_i+b_i)D_i))$ and $$st=z^{\alpha+\beta}+ \textrm{ higher order terms}.$$ Thus the subadditivity of $F[\psi_1,...,\psi_r]$ follows exactly as for $F[\psi]$ in Lemma \ref{subsub}.
\end{proof}

\begin{lemma} \label{finishedhey}
$F[\psi_1,...,\psi_r]$ is locally linearly bounded from below.  
\end{lemma}

\begin{proof}
Let $(\alpha,a)\in \Sigma(\Gamma(X))^{\circ}.$ Let $\psi_{i,p}$ be the trivializations of the metrics $\psi_i,$ then $$\sum a_i\psi_{i,p}$$ is the trivialization of $\sum a_i\psi_i.$ Let $D$ be as in the proof of Lemma \ref{lowerbound}, and choose $A$ such that $$e^{-\sum a_i\psi_{i,p}}>A.$$ Since the inequality $$e^{-\sum b_i\psi_{i,p}}>A$$ holds for all $b$ in a neighbourhood of $a,$ the lower bound follows as in the proof of Lemma \ref{lowerbound}.
\end{proof}

\begin{definition}
The global Chebyshev transform $c[\psi_1,...,\psi_r]$of the $r$-tuple $(\psi_1,...,\psi_r)$ is defined as the convex envelope of $F[\psi_1,...,\psi_r]$ on $\Sigma(\Gamma(X))^{\circ}.$
\end{definition}

\begin{proposition} \label{proorp}
For any sequence $(\alpha(k),a(k))\in \Gamma(X)$ such that $|(\alpha(k),a(k))| \to \infty$ and $$\frac{(\alpha(k),a(k))}{|(\alpha(k),a(k))|}\to (p,a)\in \Sigma(\Gamma(X))^{\circ}$$ it holds that $$\lim_{k\to \infty}\frac{F[\psi_1,...,\psi_r](\alpha(k),a(k))}{|(\alpha(k),a(k))|}=c[\psi_1,...,\psi_r](p,a).$$
\end{proposition}

\begin{proof}
By Lemma \ref{soon} and Lemma \ref{finishedhey} we can use Theorem \ref{jajatheo}, which gives us the proposition.
\end{proof}

\begin{proposition} \label{globalext}
For rational $a$, i.e $a=(a_1,...,a_r)\in \mathbb{Q}^r,$ the global Chebyshev transform $c[\psi_1,...,\psi_r](p,a)$ coincides with  $c\left[\sum a_i\psi_i\right](p),$ where the Chebyshev transform of the $\mathbb{Q}$-divisor $\sum a_iD_i$ as defined by (\ref{qdivisor}). 
\end{proposition} 

\begin{proof}
By construction it is clear that for all $(\alpha,a)\in \Gamma(X)$ we have that $$F[\psi_1,...,\psi_r](\alpha,ka)=F\left[\sum a_i\psi_i\right](\alpha,k).$$
Choose a sequence $(\alpha(k),ka)\in \Gamma(X)$ such that $$\lim_{k\to \infty}\frac{(\alpha(k),ka))}{|(\alpha(k),ka))|}=\frac{(p,a)}{|(p,a)|},$$ where we only consider those $k$ such that $ka$ is an integer. Then by Proposition \ref{proorp} we have that 
\begin{eqnarray*}
c[\psi_1,...,\psi_r](p,a)=\lim_{k\to \infty}|(p,a)|\frac{F[\psi_1,...,\psi_r](\alpha(k),ka)}{|(\alpha(k),ka)|}=\\=\lim_{k\to \infty}|(p,a)|\frac{F\left[\sum a_i\psi_i\right](\alpha(k),k)}{|(\alpha(k),ka)|}=\lim_{k\to \infty}\left(\frac{|(p,a)|k}{|(\alpha(k),ka)|}\right)c\left[\sum a_i\psi_i\right](p)=\\=c\left[\sum a_i\psi_i\right](p).
\end{eqnarray*}
\end{proof}

Since the global Okounkov body and the global Chebyshev transform are convex it follows that the formula (\ref{mainformula}) defines a continuous extension of the metric volume to the space of big metric $\mathbb{R}-$divisors.

In order to prove further regularity of the metric volume we wish to show that the formula (\ref{superformula}) still holds for the extension. First we need some preliminary lemmas.

\begin{lemma} \label{relativeenergyhom}
The function $\mathcal{E}\circ P(t\psi,t\varphi)$ is $(n+1)$-homogeneous in $t$ for $t>0,$ i.e. $$\mathcal{E}\circ P(t\psi,t\varphi)=t^{n+1}\mathcal{E}\circ P(\psi,\varphi).$$
\end{lemma}

\begin{proof}
For metrics with minimal singularities $\psi'$ and $\varphi',$ by definition of the Monge-Amp\`ere energy we have that 
\begin{eqnarray} \label{relativehom}
\mathcal{E}(t\psi,t\varphi)=\frac{1}{n+1}\int_{\Omega} (t\psi'-t\varphi')\textrm{MA}_n(t\psi',t\varphi')=\nonumber \\=\frac{t^{n+1}}{n+1}\int_{\Omega} (\psi'-\varphi')\textrm{MA}_n(\psi',\varphi')=t^{n+1}\mathcal{E}(\psi,\varphi).
\end{eqnarray}
We also observe that $t\psi'$ is a psh metric on $tL$ iff $\psi'$ is a psh metric on $L.$ Therefore we get that 
\begin{equation} \label{projhom}
P(t\psi)=tP(\psi).
\end{equation} 
Combining (\ref{relativehom}) and (\ref{projhom}) the lemma follows.
\end{proof}

\begin{lemma} \label{shitock}
Let $\psi$ and $\psi'$ be two continuous metrics on $L,$ and let $\varphi$ and $\varphi'$ be two continuous metrics on some other big line bundle $L'.$ Then the function $$\mathcal{E}\circ P(\psi+t\varphi,\psi'+t\varphi')$$ is continuous in $t$ for $t\geq 0.$ 
\end{lemma}

\begin{proof}
Since the Monge-Amp\`ere energy is homogeneous we may assume that $L'$ has a nontrivial section $S$. We let $\varphi_S$ denote the (singular) metric defined such that the function $|S|_{\varphi_S}$ is identically equal to one. Let $\psi_t$ denote the singular metric $P(\psi+t\varphi)+(1-t)\varphi_S$ and similarly let $\psi'_t$ denote the metric $P(\psi'+t\varphi')+(1-t)\varphi_S.$ Since the singular metric $\varphi_S$ is pluriharmonic outside of the zero locus it follows that 
\begin{eqnarray} \label{contintegral}
\mathcal{E}\circ P(\psi+t\varphi,\psi'+t\varphi')=\frac{1}{n+1}\int_U(P(\psi+t\varphi)-P(\psi'+t\varphi'))MA_n(\psi_t,\psi'_t),
\end{eqnarray} 
where $U$ denotes a dense Zariski open set where the metrics in question are locally bounded. 

Lemma 1.14 in \cite{Berman} tells us that the projection operator is 1-Lipschitz continuous. In our case this means that $$\sup_X |P(\psi+t\varphi)-P(\psi+t\varphi')|\leq t\sup_X |\varphi-\varphi'|.$$ Therefore $P(\psi+t\varphi)-P(\psi+t\varphi')$ is uniformly bounded on $X.$ 

For any $0\leq s\leq t$ we have that $$P(\psi+t\varphi)\geq P(\psi+s\varphi)+(t-s)P(\varphi).$$ It follows that $\psi_t$ is increasing in $t$ on the set where $P(\varphi)>\varphi_S.$ It is also easy to see that $\psi_t$ decreases to $\psi_r$ when $t$ decreases to $r,$ and that $\psi_t$ increases to $\psi_r$ a.e. when $t$ increases to $r.$ Let $U'$ denote the plurifine open set gotten by intersecting $U$ by the sets where $P(\varphi)$ and $P(\varphi')$ are greater than $\varphi_S.$ By the work of Bedford-Taylor we get that the measure 
\begin{eqnarray*}
(P(\psi+t\varphi)-P(\psi'+t\varphi'))MA_n(P(\psi+t\varphi),P(\psi'+t\varphi'))=\\=(P(\psi+t\varphi)-P(\psi'+t\varphi'))MA_n(\psi_t,\psi'_t)
\end{eqnarray*} 
restricted to $U'$ varies continuously in $t$ in the weak sense (on U'). By using $\varphi_S-C$ instead of $\varphi_S$ where $C$ is an arbitrary constant we get that the restriction of the measure to $U$ minus the zero set of $S$ varies continuously. The complement of U together with the zero set of $S$ is pluripolar, and is thus not charged by the mixed Monge-Amp\`ere measures. The total mass of $MA_n(P(\psi+t\varphi),P(\psi'+t\varphi'))$ is by \cite{Berman} equal to $n+1$ times the volume of $L+tL',$ which varies continuously with $t$ by e.g. \cite{Favre}. As in Theorem 2.6 in \cite{Berman} this implies that the integral in (\ref{contintegral}) and thus the Monge-Amp\`ere energy $\mathcal{E}\circ P(\psi+t\varphi,\psi'+t\varphi')$ is continuous in $t.$  

\end{proof}

We are now ready to prove formula (\ref{superformula}) in the setting of metrized big $\mathbb{R}$-divisors.

\begin{theorem} \label{weimain}
For big $\mathbb{R}$-divisors $\sum a_iD_i$ we have that
\begin{eqnarray} \label{realequality}
\mathcal{E}\circ P(\sum a_i\psi_i,\sum a_i\varphi_i)= \nonumber \\=n!\int_{\Delta(\sum a_iD_i)}(c[\varphi_1,...,\varphi_r](p,a)-c[\psi_1,...,\psi_r](p,a))d\lambda(p).
\end{eqnarray}
\end{theorem}

\begin{proof}
First we show that (\ref{realequality}) holds when $a\in \mathbb{Q}^r$. By the homogeneity of the Okounkov body and the Chebyshev transform we have that 
\begin{eqnarray*}n!\int_{\Delta(tL)^{\circ}}(c[t\psi]-c[t\varphi])d\lambda=t^{n+1}n!\int_{\Delta(L)^{\circ}}(c[\psi]-c[\varphi])d\lambda=\\=t^{n+1}\mathcal{E}(\varphi,\psi)=\mathcal{E}(t\varphi,t\psi),
\end{eqnarray*}
where the last equality follows from Lemma \ref{relativeenergyhom}. Then by Proposition \ref{globalext}, (\ref{realequality}) holds for $a\in \mathbb{Q}^r$. Therefore by the continuity of the Monge-Amp\`ere energy, the continuity of the global Chebyshev transform, and the fact that equation (\ref{realequality}) holds for rational $a,$ the proposition follows.
\end{proof}

 \section{Differentiability of the metric volume}

We wish to understand the behaviour of the metric volume $\textrm{vol}(L_t,\psi_t,\varphi)$ when $L_t$ and the metrics $\psi_t$ and  $\varphi_t$ vary with $t.$ In \cite{Berman} Berman-Boucksom study the case where $\psi_t$ and $\varphi_t$ are metrics on a fixed line bundle or more generally a big $\mathbb{R}$-divisor. We are interested in the case where the underlying $\mathbb{R}$-divisor $L_t$ is allowed to vary as well. As we have seen, by letting $\psi_t=\varphi_t+1$ the problem reduces to that of the variation of the volume functional. It was first proven by Boucksom-Favre-Jonsson in \cite{Favre} that the voume functional was $C^1$ on the space of big $\mathbb{R}$-divisors. In \cite{Lazarsfeld} Lazarsfeld-Musta\c{t}\u{a} reprove this differentiability result by studying the variation of the Okounkov bodies. Since our Theorem \ref{main} and Theorem \ref{weimain} states that the metric volume is given by the integration of the difference of Chebyshev transforms on the Okounkov body, we wish to use the same approach as Lazarsfeld-Musta\c{t}\u{a} did in \cite{Lazarsfeld}. The situation becomes a bit more involved, since we have to consider not only the variation of the Okounkov bodies but also the variation of the Chebyshev transforms.

To account for the variation of the Chebyshev transform when the underlying line bundle changes it becomes necessary to consider not only continuous metrics but also metrics with certain singularities. 

Let $S$ denote a section of an ample line bundle $A.$ As above we let $\varphi_S$ denote the (singular) metric defined such that the function $|S|_{\varphi_S}$ is identically equal to one. Let also $\Psi$ be some fixed continuous positive metric on $A.$ For any number $R$ we denote by $\varphi_{S,R}$ the metric $$\varphi_{S,R}:=\max(\varphi_S,\Psi-R).$$ 

\begin{lemma} \label{prrr}
Let $\psi$ be a continuous metric on a big line bundle $L,$ and let $t>0$ be such that $L-tA$ is still big. For $R \gg 0$ we have that $$P(\psi-t\varphi_{S,R})=P(\psi-t\varphi_S).$$
\end{lemma}

\begin{proof}
That $$P(\psi-t\varphi_{S,R})\leq P(\psi-t\varphi_S)$$ is clear since $$\psi-t\varphi_{S,R}\leq \psi-t\varphi_S.$$ $P(\psi-t\varphi_S)$ is psh, therefore upper semicontinuous by definition, which means that it is locally bounded from above. Thus locally we can find $R \gg 0$ such that $$\psi-t(\Psi-R)\geq P(\psi-t\varphi_S).$$ But we have assumed that our manifold $X$ is compact, so there exists an $R$ such that $\psi-t(\Psi-R)$ dominates $P(\psi-t\varphi_S)$ on the whole of $X.$ The same must be true for $\psi-t\varphi_{S,R}.$ By definition $P(\psi-t\varphi_{S,R})$ dominates all psh metrics less or equal to $\psi-t\varphi_{S,R},$ in particular it must dominate $P(\psi-t\varphi_S).$
\end{proof}

\begin{lemma} \label{singular}
If $L$ is integral, i.e. a line bundle, then for large enough $R$ the function $F[\psi-t\varphi_{S,R}]$ is independent of $R,$ and we will use $F[\psi-t\varphi_S]$ to denote this function.
\end{lemma}

\begin{proof}
This follows the fact that for all metrics $\varphi$ and all sections $s$ it holds that $$\sup_{x\in X}\{|s(x)|^2e^{-\varphi(x)}\}=\sup_{x\in X}\{|s(x)|^2e^{-P(\varphi)(x)}\},$$ see e.g. \cite{Berman}. 
\end{proof}

From Lemma \ref{singular} it follows that the Chebyshev transform $c[\psi-t\varphi_S]$ is well-defined, also for $\mathbb{R}-$divisors, and that Proposition \ref{propidopp} holds in this case. The formula for the Monge-Amp\`ere energy as the integral of Chebyshev transforms will also still hold.

\begin{proposition} 
For any continuous metric $\varphi$ on $L-tA$ it holds that
\begin{eqnarray} \label{tratralepp}
\mathcal{E}\circ P(\psi-t\varphi_S,\varphi)=\\=n!\int_{\Delta(L-tA)^{\circ}}c[\varphi]-c[\psi-t\varphi_S]d\lambda.
\end{eqnarray}
\end{proposition}

\begin{proof}
For integral $L,$ choose an $R\gg 0$ such that $$P(\psi-t\varphi_{S,R})=P(\psi-t\varphi_S).$$ Then (\ref{tratralepp}) follows in this case from Theorem \ref{main} and Lemma \ref{singular}. By homogeneity (\ref{tratralepp}) holds for rational $L,$ and by continuity for arbitrary big $\mathbb{R}$-divisors.
\end{proof}

Theorem B in \cite{Berman} states that the Monge-Amp\`ere energy is differentiable when the metrics correspond to a fixed big line bundle. By the comment in the beginning of section 4 in \cite{Berman} this holds more generally for big $(1,1)$ cohomology classes, e.g. $\mathbb{R}$-divisors. We thus have the following.

\begin{theorem} \label{Berman}
Let $\psi$ and $\varphi$ be continuous metrics on a big $\mathbb{R}$-divisor $D$, and let $u$ be a continuous function. Then the function $$f(t):=\textrm{vol}(D,\psi+tu,\varphi)$$ is differentiable, and $$f'(0)=\int_{\Omega}u\textrm{MA}(P(\psi_0)).$$
\end{theorem}

We also need to consider the case where $$\psi_t=\psi_0+t(\Phi-\varphi_S),$$ where $\Phi$ is some continuous metric on $A.$ 

We state and prove a slight variation of Lemma 3.1 in \cite{Witt}.
 
 \begin{lemma} \label{witt}
 Let $f_k$ be a sequence of concave functions on the unit interval increasing pointwise to a concave function $g.$ Then $$g'(0)\leq \liminf_{k\to \infty}f'_k(0),$$ allowing the possibility that $f'_k(0)$ and $g'(0)$ are plus infinity. 
 \end{lemma}

\begin{proof}
Since $f_k$ is concave we have that $$f_k(0)+f'_k(0)t\geq f_k(t)$$ hence $$\liminf_{k\to \infty}tf'_k(0)\geq g(t)-g(0).$$ The lemma follows by letting $t$ tend to zero.
\end{proof}

\begin{lemma} \label{wiwi}
The function $$f(t):=\mathcal{E}\circ P(\psi_0+t(\Phi-\varphi_S),\varphi)$$ is concave for $t\geq 0$ and for $t>0$ we have that
\begin{equation} \label{ineqder}
\frac{d}{dt}_{|t+}f\leq \int_{\Omega}(\Phi-\varphi_S)\textrm{MA}(P(\psi_0+t(\Phi-\varphi_S)))\leq \frac{d}{dt}_{|t-}f.
\end{equation}
\end{lemma}

\begin{proof} 
Without loss of generality we can assume that $\Phi-\varphi_S>0.$ Thus $P(\psi_0+t(\Phi-\varphi_S))$ and therefore $f(t)$ is increasing in $t.$ 

Let us denote $\Phi-\varphi_S$ by $u,$ and let $$u_k:=\Phi-\varphi_{S,k}.$$ Let $f_k$ denote the function $$f_k(t):=\mathcal{E}\circ P(\psi_0+tu_k,\varphi).$$ By e.g. \cite{Berman} the functions $f_k$ are concave, and by Theorem \ref{Berman} they are differentiable with $$f'_k(t)=\int_{\Omega}u_k \textrm{MA}(P(\psi_0+tu_k)).$$ Clearly $f_k$ is increasing in $k$ and by Lemma \ref{prrr} $f_k$ increases pointwise to $f.$ It follows that $f$ is concave. 

By Lemma \ref{prrr} and monotone convergence we get that $f'_k(t)$ converges to $$\int_{\Omega}u\textrm{MA}(P(\psi_0+tu))$$ pointwise. The inequalities (\ref{ineqder}) follow from applying Lemma \ref{witt} to $f$ and its reflection. 

\end{proof}

We will also need an integration by parts formula involving $\varphi_S,$ which generalizes Proposition 4.7 in \cite{Berman}.

\begin{lemma} \label{partialint}
Let $\varphi$ and $\varphi'$ be continuous metrics on a big $\mathbb{R}$-divisor $L.$ Let $\psi$ be a continuous psh metric on an ample line bundle $A,$ and let $S\in H^0(A)$ be a section such that its zero set variety $Y$ is a smooth submanifold not contained in the augmented base locus of $L$. Then it holds that 
\begin{eqnarray*}
\int_X (\psi-\varphi_S)(\textrm{MA}(P(\varphi))-\textrm{MA}(P(\varphi')))=\\=\int_X (P(\varphi)-P(\varphi'))dd^c \psi \wedge \textrm{MA}_{n-1}(P(\varphi),P(\varphi'))-n\mathcal{E}_Y(P(\varphi)_{|Y},P(\varphi')_{|Y}).
\end{eqnarray*}  
\end{lemma}

\begin{proof}
Following the proof of Proposition 4.7 \cite{Berman} we observe that $$\textrm{MA}(P(\varphi))-\textrm{MA}(P(\varphi'))=dd^c(P(\varphi)-P(\varphi'))\wedge \textrm{MA}_{n-1}(P(\varphi),P(\varphi')).$$ 
The lemma will follow by the Lelong-Poincar\'e formula as soon as we establish that 
\begin{eqnarray*}
\int_X (\psi-\varphi_S)dd^c (P(\varphi)-P(\varphi'))\wedge \textrm{MA}_{n-1}(P(\varphi),P(\varphi'))=\\=\int_X (P(\varphi)-P(\varphi'))dd^c (\psi-\varphi_S)\wedge \textrm{MA}_{n-1}(P(\varphi),P(\varphi')),
\end{eqnarray*}
which is an integration by parts formula. By \cite{BEGZ} we may integrate by parts when the functions are differences of quasi-psh metrics with minimal singularities. We denote by $u_k$ the quasi-psh metric with minimal singularities $\psi-\varphi_{S,k}$ and get that
\begin{eqnarray} \label{impequation}
\int_X u_kdd^c (P(\varphi)-P(\varphi'))\wedge \textrm{MA}_{n-1}(P(\varphi),P(\varphi'))=\\=\int_X (P(\varphi)-P(\varphi'))dd^c u_k\wedge \textrm{MA}_{n-1}(P(\varphi),P(\varphi')).
\end{eqnarray}
Let $U$ be the dense Zariski open set where $P(\varphi)$ and $P(\varphi')$ is locally bounded. As in the proof of Lemma \ref{shitock} we get that $dd^c u_k\wedge \textrm{MA}_{n-1}(P(\varphi),P(\varphi'))$ converge weakly to $$(P(\varphi)-P(\varphi')dd^c\varphi_S\wedge \textrm{MA}_{n-1}(P(\varphi),P(\varphi'))$$ on $U.$ The integral $$\int_U dd^c u_k\wedge \textrm{MA}_{n-1}(P(\varphi),P(\varphi'))$$ equal to $n+1$ times the restricted volume $\langle L^{n-1}A\rangle,$ see \cite{Favre} and \cite{BEGZ} and thus independent of $k.$ It is easily seen that also $$\int_U dd^c \varphi_S\wedge \textrm{MA}_{n-1}(P(\varphi),P(\varphi'))=(n+1)\langle L^{n-1}A\rangle,$$ see e.g. \cite{Hisamoto}. Since the restricted volume varies continuously with $L$ (see \cite{Favre}) it follows as in the proof of Lemma \ref{shitock} that the right hand side of equation (\ref{impequation}) converges to $$\int_Y(P(\varphi)-P(\varphi'))\wedge \textrm{MA}_{n-1}(P(\varphi),P(\varphi'))=\mathcal{E}_Y(P(\varphi)_{|Y},P(\varphi')_{|Y})$$ when $k$ tends to infinity.

Clearly by monotone convergence $$\int_X u_k\textrm{MA}(P(\varphi))$$ converges to 
\begin{equation} \label{diffeq11}
\int_X (\psi-\varphi_S)(\textrm{MA}(P(\varphi)),
\end{equation} 
but we need to show that it is finite to conclude that $$\int_X u_k(\textrm{MA}(P(\varphi))-\textrm{MA}(P(\varphi'))$$ converges to 
\begin{equation} \label{diffeq21}
\int_X (\psi-\varphi_S)(\textrm{MA}(P(\varphi))-\textrm{MA}(P(\varphi')).
\end{equation}

From the inequality (\ref{ineqder}) we see that $$\int_X (\psi-\varphi_S)(\textrm{MA}(P(\varphi_0-r\varphi_S)$$ is finite when $\varphi_0$ is a continuous metric on $L+rA.$ Lettin $\varphi$ be $\varphi_0-r\varphi_{S,k}$ where $k$ is chosen so that $P(\varphi_0-r\varphi_{S,k})=P(\varphi_0-r\varphi_S)$ we see that $$\int_X (\psi-\varphi_S)(\textrm{MA}(P(\varphi))$$ is finite for at least one continuous metric $\varphi.$ But using equation (\ref{impequation}) we see that the absolute value of the difference (\ref{diffeq21}) is bounded by a uniform constant times the supremum of $|P(\varphi)-P(\varphi')|,$ which is bounded. We conclude that (\ref{diffeq11}) is bounded for all continuous metrics $\varphi,$ and therefore the Lemma follows from applying monotone convergence.
 
\end{proof}

\begin{corollary} \label{importcor}
The function $$f(t):=\mathcal{E}\circ P(\psi_0+t(\Phi-\varphi_S),\varphi)$$ is continuosly differentiable with $$f'(t)=\int_{\Omega}(\Phi-\varphi_S)\textrm{MA}(P(\psi_0+t(\Phi-\varphi_S)),$$ where for $t=0,$ $f'(0)$ here denotes the right derivative.
\end{corollary}

\begin{proof}
Using the integration by parts formula one argues as in the proof of Lemma \ref{shitock} and concludes that $$\int_{\Omega}(\Phi-\varphi_S)\textrm{MA}(P(\psi_0+t(\Phi-\varphi_S))$$ varies continuosly with $t.$ Since $f'(t)$ is decreasing, the corollary follows from the inequality (\ref{ineqder}). 
\end{proof}

Assume that we have chosen our coordinates $z_1,...,z_n$ centered at $p$ such that $$z_1=0$$ is a local equation for an irreducible variety $Y$ not contained in the augmented base locus of $L.$ Assume also that $Y$ is the zero-set of a holomorphic section $S\in H^0(A)$ of an ample line bundle $A.$ Then by Theorem 4.24 in \cite{Lazarsfeld} the Okounkov bodies of $L$ and $L+tA$ with respect to these coordinates are related in the following way $$\Delta(L)=(\Delta(L+tA)-te_1)\cap (\mathbb{R}_+)^n.$$ There is also correspondence between the Chebyshev transforms of metrics on $L$ and $L+tA.$

\begin{proposition} \label{trollalla}
Let $A$ and $S$ be as above. Suppose also that we have chosen the holomorphic coordinates so that $z_1=S$ locally. Then for $a> r$ it holds that 
\begin{eqnarray} \label{sososo}
c_L[\psi](a,\alpha)-c_L[\varphi](a,\alpha)=\nonumber \\ =c_{L-rA}[\psi-r\varphi_S](a-r,\alpha)-c_{L-rA}[\varphi - r\varphi_S](a-r,\alpha).
\end{eqnarray}
\end{proposition}

\begin{proof}
First assume that $L$ is integral. Since we have that locally $S=z_1,$ for $t \in H^0(kL),$ $$t=z^{k(a,\alpha)}+\textrm{higher order terms},$$ if and only if $$\frac{t}{S^{rk}}=z^{k(a-r,\alpha)}+\textrm{higher order terms}.$$ We also have that $$\sup_{x\in X} \{|t(x)|^2 e^{-k\varphi(x)}\} =\sup_{x\in X} \{\frac{|t(x)|^2}{|s^{rk}(x)|^2}e^{-k(\varphi(x)-r\ln |s(x)|^2)}\}.$$ Thus (\ref{sososo}) holds for integral $L$. By the homogeneity and continuity of the Chebyshev transform it will therefore hold for big $\mathbb{R}$-divisors.
\end{proof}

We are now ready to state and prove our differentiability theorem for the metric volume. 

\begin{theorem}
Let $L_i,$ $i=1,...,m$ be a collection of line bundles, and for each $i$ let $\psi_i$ and $\varphi_i$ be two continuous metrics on $L_i.$ Denote $\sum a_iL_i$ by $L_a,$ $\sum a_i\psi_i$ by $\psi_a$ and $\sum a_i\varphi_i$ by $\varphi_a.$ Let $O$ denote the open cone in $\mathbb{R}^m$ such that $a\in O$ iff $L_a$ is big. Then the function $$f(a):=\textrm{vol}(L_a,\psi_a,\varphi_a)$$ is $\mathcal{C}^1$ on $O.$
\end{theorem}

\begin{proof}
Let $a$ be a point in $O,$ and let $L=L_a.$ Denote $\psi_a$ by $\psi$ and $\varphi_a$ by $\varphi.$  Let us consider the (possible) partial derivative of $f$ at $a$ in the $x_1$-direction. Since any line bundle can be written as the difference of very ample line bundles, without loss of generality we can assume that $L_1=A$ is ample and has a section $S$ defining a smooth hypersurface $Y$ not contained in the augmented base locus of $L$. Let us denote the metrics on $A$ by $\Psi$ and $\Phi$ in order to avoid confusion. We consider the restricted function $$f(t):=\textrm{vol}(L+tA,\psi + t\Psi, \varphi+t\Phi).$$ We claim that $f$ is differentiable at $t=0$, and that the derivative varies continuously with $L,$ $\psi$ and $\varphi.$ This will imply that the function $f$ was $\mathcal{C}^1$ on the whole of $O.$

We choose local holomorphic coordinates such that $z_1=S.$ Recall that the Okounkov bodies of $L$ and $L+tA$ are related in the following way 
\begin{equation} \label{snartklart}
\Delta(L)=(\Delta(L+tA)-te_1)\cap (\mathbb{R}_+)^n.
\end{equation}

Let $\Delta(L)_r$ denote the fiber over $r$ of the projection of the Okounkov body down to the first coordinate, i.e. $$\Delta(L)_r:=\Delta(L)\cap (\{r\}\times \mathbb{R}^{n-1}).$$ Then one may write equation (\ref{snartklart}) as 
\begin{equation} \label{snartsnart}
\Delta(L+tA)=\cup_{0\leq r\leq t}\Delta(L+tA)_r\cup(\Delta(L)+te_1).
\end{equation}

Furthermore the metric volume is given by integration of the Chebyshev transforms over the Okounkov bodies. Using (\ref{snartsnart}) and Proposition \ref{trollalla} we get that

\begin{eqnarray*}
\textrm{vol}(L+tA,\psi + t\Psi, \varphi+t\Phi)=\\=n!\int_{\Delta(L+tA)^{\circ}}c[\varphi+t\Phi]-c[\psi+t\Psi]d\lambda=\\ =n!\int_{r=0}^t \int_{\Delta(L+tA)^{\circ}_r}c[\varphi+t\Phi](r,\alpha)-c[\psi+t\Psi](r,\alpha)d\alpha dr + \\ +n!\int_{\Delta(L)^{\circ}}c[\varphi+ t(\Phi -\varphi_S)]-c[\psi+t(\Psi-\varphi_S)]dp=\\ =n!\int_{r=0}^t \int_{\Delta(L+tA)^{\circ}_r}c[\varphi+t\Phi](r,\alpha)-c[\psi+t\Psi](r,\alpha)d\alpha dr+\\+ \mathcal{E}_L\circ P(\psi + t(\Psi -\varphi_S), \varphi + t(\Phi -\varphi_S)).
\end{eqnarray*}

As in Corollary \ref{corzero} the global Chebyshev transforms $c[\psi,\Psi]$ and $c[\varphi,\Phi]$ will have continuous extensions to the interior  of the zero fiber of the corresponding global Okounkov body, which simply consists of the zero fibers of $\Delta(L+tA).$ Hence by the fundamental theorem of calculus and Corollary \ref{importcor} it follows that $f$ is right-differentiable.

We get that 
\begin{eqnarray*}
\frac{d}{dt}_{|_{0+}} \textrm{vol}(L+tA,\psi + t\Psi, \varphi+t\Phi)=\\= n!\int_{\Delta(L)^{\circ}_0}c[\varphi](0,\alpha)-c[\psi](0,\alpha)d\alpha+\\ +\frac{d}{dt}_{|_{0+}} \mathcal{E}_{L}\circ P(\psi + t(\Psi-\varphi_S), \varphi+t(\Phi-\varphi_S)).
\end{eqnarray*}

The first term depends continuously on the data since $\Delta(L)^{\circ}_0$ depends continuously on $L$ and since the global Chebyshev transforms and its extensions are continuous.

Let us look at the second term. Because of the cocycle property of the Monge-Amp\`ere energy, we only need to consider two cases, one where $\psi=\varphi,$ and the other one where we let $\psi \neq \varphi$ but instead assume that $\Psi=\Phi.$

First assume that $\psi=\varphi$. We get that
\begin{eqnarray} \label{hurrahso}
\frac{d}{dt}_{|_0+} \mathcal{E}_{L}\circ P(\psi + t(\Psi-\varphi_S2),\psi+t(\Phi-\varphi_S))=\nonumber \\=\int_X (\Psi-\varphi_S)\textrm{MA}(P(\psi))-\int_X (\Phi-\varphi_S)\textrm{MA}(P(\psi))=\nonumber \\=\int_X (\Psi-\Phi)\textrm{MA}(P(\psi)).
\end{eqnarray}

As in Lemma \ref{shitock} this will depend continuously on $\psi.$

Now let $\psi \neq \varphi$ but instead assume that $\Psi=\Phi'$ is some metric on $A.$ By Lemma \ref{wiwi}, the cocycle property, and the integration by parts formula in Lemma \ref{partialint} we have that

\begin{eqnarray} \label{troski}
\frac{d}{dt}_{|_0+} \mathcal{E}_{L}\circ P(\psi + t(\Psi-\varphi_S), \varphi+t(\Psi-\varphi_S))=\nonumber \\=\int_X (\Psi-\varphi_S)(\textrm{MA}(P(\psi))-\textrm{MA}(P(\varphi)))=\\=\int_X (P(\psi)-P(\varphi))dd^c \psi \wedge \textrm{MA}_{n-1}(P(\psi),P(\varphi))-n\mathcal{E}_Y(P(\psi)_{|Y},P(\varphi)_{|Y}).
\end{eqnarray} 

Arguing as in the proof of Lemma \ref{shitock} we get that both these terms depend continuously on the data.

By Lemma \ref{shitock} the function $f(t)$ is continuous and we have seen that it is also continuosly right differentiable. An elementary application of the mean value for right differentiable functions yields that any continuous and continuously right differentiable function is in fact differentiable. Thus $f(t)$ is continuously differentiable for $t>0$ and since the choice of $L$ was arbitrary it is differentiable in a neighbourhood of zero as well. Since the derivative depended continuosly on the data, and thus on the point $a\in O$ it follows that $f$ is $\mathcal{C}^1$ in $O.$
\end{proof}

$ $\\
David Witt Nystr\"{o}m\\
University of Cambridge\\
Department of Pure Mathematics and Mathematical Statistics \\
 CB3 0WB Cambridge\\
United Kingdom \\
email: D.WittNystrom@dpmms.cam.ac.uk or danspolitik@gmail.com

\end{document}